\documentclass[11pt,reqno]{amsart}
\usepackage{amssymb}
\usepackage{amsmath, amssymb, amsthm}
\usepackage{mathrsfs}

\usepackage{color}

\newtheorem{theorem}{Theorem}[section]
\newtheorem{definition}{Definition}[section]
\newtheorem{lemma}{Lemma}[section]
\newtheorem{remark}{Remark}[section]

\newtheorem{corollary}{Corollary}[section]
\numberwithin{equation}{section}

\begin{document}
\title[Shallow water equations]{On classical solutions to 2D Shallow water equations with degenerate viscosities}

\author{yachun li  }
\address[Y. C. Li]{Department of Mathematics \& Key Lab of Scientific Engineering Computing (MOE),  Shanghai Jiao Tong University,
Shanghai 200240, P.R.China} \email{\tt ycli@sjtu.edu.cn}
\author{ronghua pan  }
\address[R. H. Pan]{School of Mathematics, Georgia Institute of Technology,
Atlanta 30332, USA} \email{\tt panrh@math.gatech.edu }

\author{shengguo Zhu}
\address[S. G. Zhu]{Department of Mathematics, Shanghai Jiao Tong University,
Shanghai 200240,  P. R. China;
School of Mathematics, Georgia Institute of Technology, Atlanta 30332, USA}
\email{\tt zhushengguo@sjtu.edu.cn}

\begin{abstract}
2D shallow water equations have degenerate viscosities proportional to surface height, which vanishes in many physical considerations, say, when the initial total mass, or energy are finite. Such a  degeneracy is a highly challenging
obstacle for development of well-posedness theory, even local-in-time theory remains open for long time.
In this paper, we will address this open problem  with some new perspectives, independent of the celebrated
BD-entropy  \cite{bd2, bd}. After exploring some interesting structures of most models of 2D shallow water equations, we introduced
a proper notion of solution class, called regular solutions, and identified a class of initial data with finite
total mass and energy, and established the local-in-time well-posedness of this class of smooth solutions.
The theory is applicable to most relatively physical shallow water models, broader than those with BD-entropy structures. Later, a Beale-Kato-Majda type blow-up criterion is also established. This paper is mainly based on our early preprint \cite{LPZ}.
\end{abstract}

\date{March 08, 2015}

\maketitle

\section{Introduction}\ \\

In this paper, we aim at a long standing open problem in the theory of shallow water equations, that is to prove the local-in-time well-posedness of classical solutions with finite total mass and /or energy. For this purpose, we consider the following Cauchy problem of shallow water equations:
\begin{equation}
\label{eq:1.1ww}
\begin{cases}
h_t+\text{div}(h U)=0,\\[6pt]
(h U)_t+\text{div}(h U\otimes U)
  +\nabla h^2 =\mathcal{V}(h, U),\\[6pt]
(h,U)|_{t=0}=(h_0(x),U_0(x)),\quad x\in \mathbb{R}^2,\\[6pt]
(h,U)\rightarrow (0,0) \quad \text{as } \quad |x|\rightarrow \infty,\quad t> 0.
\end{cases}
\end{equation}
Here $h$ denotes the height of the free surface; $U=(U^{(1)},U^{(2)})^\top\in \mathbb{R}^2$ is the mean horizontal velocity of fluids. $\mathcal{V}(h, U)$ is the viscous term. There are several different viscous terms imposed, such as $\text{div}(h D(U)), \ \text{div}(h\nabla U), \ h\Delta U, \ \Delta (hU)$ (see \cite{bd} and \cite{lions}), where $D(U)=\frac12 [\nabla U+\nabla^T U]$. In particular, the case for $\mathcal{V}=\text{div}(h\nabla U)$ is corresponding to the well-known viscous Saint-Venant model \cite{gp}.  The derivation of Gent \cite{Gent} suggests $\mathcal{V}=\text{div}(h D(U))$. A more recent careful derivation by Marche \cite{Mar} and by Bresch-Noble \cite{BN1, BN2} suggests that
\begin{equation}\label{newvis}
\mathcal{V}(h, U)=\text{div}(2h D(U)+2h \text{div} U \mathbb{I}_2).
\end{equation}

When one replaces $h$ with $\rho$, the density of a compressible fluid, the shallow water equations
\eqref{eq:1.1ww} turns to isentropic Navier-Stokes equations with density dependent
viscosities. In order to avoid the choices amongst many interesting shallow water models, and
 to establish a solid local well-posedness theory of classical solutions to \eqref{eq:1.1ww} with broad spectrum of applications to most of the physically relevant shallow water models, as well as for possible applications to  polytropic fluids with density dependent viscosities (\cite{chap, tlt}), we will mainly work on the following
compressible Navier-Stokes equations  for the isentropic flow when viscosity coefficients, shear and bulk, are degenerate, i.e.,
\begin{equation}
\label{eq:1.1}
\begin{cases}
\rho_t+\text{div}(\rho u)=0,\   x\in \mathbb{R}^2, \\[4pt]
(\rho u)_t+\text{div}(\rho u\otimes u)
  +\nabla
   P =\text{div} \mathbb{T}\\
(\rho,u)|_{t=0}=(\rho_0(x),u_0(x)),\quad x\in \mathbb{R}^2,\\
(\rho,u)\rightarrow (0,0) \quad \text{as } \quad |x|\rightarrow +\infty,\quad t> 0.
\end{cases}
\end{equation}
 $\rho$ is the density, and $u=(u^{(1)},u^{(2)})^\top\in \mathbb{R}^2$ is the velocity of the fluid. We assume that the pressure  $P$ satisfies
\begin{equation}
\label{eq:1.2}
P=A\rho^{\gamma}, \quad \gamma> 1,
\end{equation}
where $A$ is a positive constant, $\gamma$ is the adiabatic exponent. $\mathbb{T}$ denotes the viscosity stress tensor with the following form
\begin{equation}
\label{eq:1.3}
\mathbb{T}=\mu(\rho)(\nabla u+(\nabla u)^\top)+\lambda(\rho) \text{div}u\mathbb{I}_2.
\end{equation}
 Here $\mathbb{I}_2$ is the $2\times 2$ identity matrix, $\mu(\rho)=\alpha\rho$ is the shear viscosity coefficient, $\lambda(\rho)=\beta\rho$, and $\mu(\rho)+\lambda(\rho)$ is the bulk viscosity coefficient, $\alpha$ and $\beta$ are both constants satisfying
 \begin{equation}\label{10000}\alpha>0,\quad  \alpha+\beta\geq 0.
 \end{equation}
Clearly, when $\gamma=2$, \eqref{eq:1.1}--\eqref{10000} cover most of shallow water models mentioned
above. We will establish the corresponding results on \eqref{eq:1.1}--\eqref{10000}, then reveal the applications to various shallow water models.

System \eqref{eq:1.1}  has received extensive attention recently. For some interesting development, we refer
the readers to  \cite{bd2}, \cite{bd}, \cite{sz2}, \cite{hailiang}, \cite{taiping}, \cite{tyc2}, \cite{zyj}. We remark that, in spite of these significant achievements, a lot of questions remain open, including the local well-posedness of classical solutions in multiple dimensions when $\inf \rho_0=0$. Our
results in this paper is a very first step toward this direction.

 In the development of theories on shallow water equation and/or \eqref{eq:1.1}, there is a remarkable
discovery of a new mathematical entropy function by Bresch and Desjardins \cite{bd3} for
 $\lambda(\rho)$ and $\mu(\rho)$ satisfying $\lambda(\rho)=2(\mu'(\rho)\rho-\mu(\rho))$. This new entropy offers a nice estimate
$$
\mu'(\rho)\frac{\nabla \rho}{\sqrt{\rho}}\in L^\infty ([0,T];L^2(\mathbb{R}^2))
$$
provided that $ \mu'(\rho_0)\frac{\nabla \rho_0}{\sqrt{\rho_0}}\in L^2(\mathbb{R}^2)$.  For our problem, such an entropy exists when $\beta=0$. However, BD-entropy structure ruled out many interesting models, including
the one in \eqref{newvis}. In this paper, we will introduce a new approach toward regularity estimates without using BD-entropy. This helps the theory to be applicable to a much broader class of models.

When density function connects to vacuum continuously locally or in the far field, the classical frameworks for constructing local classical solutions (c.f. \cite{nash}) for the case when $\rho_0$ has a positive lower bound are not applicable. In fact, there are two major challenges. The first one lies at the degeneracy of time evolution in momentum equations. We notice that the leading coefficient of  $u_{t}$ in momentum equations vanishes at vacuum, and this leads to infinitely many ways to define velocity (if it exits) when
vacuum appears.  Mathematically, this degeneracy leads to an essential difficulty in determining the velocity when
vacuum occurs, since it is difficult to find a reasonable way to extend the definition of velocity into vacuum region.  Physically, it is not clear how to define the fluid velocity when there is no fluid at vacuum.
When viscosity coefficients are constants, a remedy was suggested in a series of papers by Cho et. al. (see \cite{CK3}, \cite{CK}, \cite{guahu}), where they imposed  initially a {\it compatibility condition}
\begin{equation*}
-\text{div} \mathbb{T}_{0}+\nabla P(\rho_0)=\sqrt{\rho_{0}} g
\end{equation*}
for some $g\in L^2$. Under this initial layer compatibility condition, a local theory was established successfully; see also \cite{luoluo}. Such a local solution was further extended globally by Huang-Li-Xin \cite{HX1}, when initial energy is small. Some similar results can also be found in \cite{HX2} and \cite{erwei}.  For the case of shallow water models, and \eqref{eq:1.1}, another essential difficulty arises due to vacuum in our case. We note that under assumption \eqref{eq:1.3},  viscosity coefficients vanish as density function connects to vacuum continuously in the far field.  This degeneracy gives rise to some difficulties in our analysis because of the less regularizing effect of the viscosity on the solutions.  This is one of the major obstacles preventing us from utilizing a similar remedy proposed by Cho et. al. for the case of constant viscosity coefficients. In order to overcome these difficulties, we observe that, assumping $\rho>0$, the momentum equations can be rewritten as
\begin{equation}
u_t+u\cdot\nabla u +\frac{2A\gamma}{\gamma-1}\rho^{\frac{\gamma-1}{2}}\nabla\rho^{\frac{\gamma-1}{2}}+Lu=\Big(\frac{\nabla \rho}{\rho}\Big) \cdot Q(u),
\end{equation}
where the so-called Lam\'e operator $ L$ and  operator $Q $ are
\begin{equation}
Lu=-\alpha\triangle u-(\alpha+\beta)\nabla \mathtt{div}u,\ Q(u)=\alpha(\nabla u+(\nabla u)^\top)+\beta\mathtt{div}u\mathbb{I}_2.
\end{equation}
Therefore, the two quantities $\rho^{\frac{\gamma-1}{2}}$ (local sound speed when multiplied by constant $\sqrt{A\gamma}$) and $\frac{\nabla \rho}{\rho}$ play significant roles in our analysis on velocity. With the help of this observation, we introduce a proper class of solutions, called regular solutions, for our problem.

\begin{definition}[\text{\textbf{Regular solution to Cauchy problem (\ref{eq:1.1})}}]\label{d1}\
 Let $T> 0$ be a finite constant. A solution $(\rho,u)$ to Cauchy problem (\ref{eq:1.1}) is called a regular solution in $ [0,T]\times \mathbb{R}^2$ if $(\rho,u)$ satisfies
\begin{equation*}\begin{split}
&(\textrm{A})\quad \rho> 0,\ \rho\in C^1([0,T]\times \mathbb{R}^2),\   \rho^{\frac{\gamma-1}{2}}\in C([0,T]; H^3), \ (\rho^{\frac{\gamma-1}{2}})_t \in C([0,T]; H^2); \\
&(\textrm{B})\quad  \nabla \rho/\rho\in C([0,T] ; L^6 \cap  D^1\cap D^2),\ (\nabla \rho/\rho)_t \in C([0,T]; H^1);\\
&(\textrm{C})\quad u\in C([0,T]; H^3)\cap L^2([0,T] ; H^4), \ u_t \in C([0,T]; H^1)\cap L^2([0,T] ; D^2);\\
&(\textrm{D})\quad u_t+u\cdot\nabla u +Lu=(\nabla \rho/\rho) \cdot Q(u),  \mathtt{for}\ |x|\to +\infty.
\end{split}
\end{equation*}
\end{definition}

Here and throughout this paper, we adopt the following simplified notations for the standard homogeneous and inhomogeneous Sobolev spaces:
\begin{equation*}\begin{split}
&D^{k,r}=\{f\in L^1_{loc}(\mathbb{R}^2): |f|_{D^{k,r}}=|\nabla^kf|_{L^r}<+\infty\}, \ D^k=D^{k,2} ,\\[6pt]
 & |f|_{D^{k,r}}=\|f\|_{D^{k,r}(\mathbb{R}^2)}, \ \|f\|_s=\|f\|_{H^s(\mathbb{R}^2)},\quad |f|_p=\|f\|_{L^p(\mathbb{R}^2)}.
\end{split}
\end{equation*}
 A detailed study of homogeneous Sobolev spaces can be found in \cite{gandi}.

\begin{remark}\label{r1}
From Definition \ref{d1} we know  that
$$\rho^{\frac{\gamma-1}{2}}\in C([0,T]; H^3),\  \nabla \rho/\rho\in C([0,T] ;L^6\cap D^1\cap D^2),$$
which imply  that $\nabla \rho/\rho\in L^\infty$. Therefore, our regular solution does not contain local vacuum, but vacuum occurs in the far field.
\end{remark}

Although it is very different, this notion of regular solutons is motivated by Makino, Ukai and Kawashima \cite{tms1} where the existence of local classical solutions to compressible Euler equations was established, using an interesting symmetrizer and the classical theory for symmetric hyperbolic systems. Similar notion was used also in
\cite{sz1},  \cite{tpy}, \cite{makio}, \cite{tyc2} and \cite{sz3}. Here, facing our much harder degenerate
hyperbolic-parabolic systems \eqref{eq:1.1}, this specialized notion of regular solutions selects velocity in a physically reasonable way when density approaches to vacuum. With the help of this notion of solutions, the  momentum equations can be reformulated into a quasi-linear parabolic system with some special source terms. Thus the problem becomes trackable through a successful linearization and  approximation process.

 In this paper, we impose the following regularity conditions on the initial data:
\begin{equation}\label{th78}
\rho_0> 0,\quad (\rho^{\frac{\gamma-1}{2}}_0, u_0)\in H^3, \quad \nabla \rho_0/\rho_0\in L^6 \cap  D^1\cap D^2.
\end{equation}

\begin{remark}\label{r2} We remark that \eqref{th78} identifies a class of admissible initial data that
provides unique solvability to our problem \eqref{eq:1.1}. This set
of initial data contains a large class of functions, for example,
$$
\rho_0(x)=\frac{1}{1+|x|^{2\sigma}}, \quad u_0(x)\in C_0^3(\mathbb{R}^2), \quad x\in \mathbb{R}^2,
$$
where $\sigma>\max\{1, \frac{1}{\gamma-1}\}$. On the other hand, the restrictions in \eqref{th78} also
reflect the difficulty of the problem, which is not completely settled even with our efforts in this paper.
\end{remark}

Now we are ready to state our main results.

\begin{theorem}\label{th2} If the initial data $( \rho_0, u_0)$ satisfy the regularity conditions in \eqref{th78}, then there exists a time $T_*>0$ and a unique regular solution $(\rho, u)$ to the Cauchy problem (\ref{eq:1.1}), satisfying
\begin{equation}\label{reg11}\begin{split}
& \rho^{\frac{\gamma-1}{2}} \in C([0,T_*];H^3),\ (\rho^{\frac{\gamma-1}{2}})_t \in C([0,T_*];H^2) ,\\
&\nabla \rho/\rho \in C([0,T_*]; L^6\cap D^1\cap D^2),\ (\nabla \rho/\rho)_t \in C([0,T_*]; H^1),\\
& u\in C([0,T_*]; H^3)\cap L^2([0,T_*] ; H^4), \ u_t \in C([0,T_*]; H^1)\cap L^2([0,T_*] ; D^2),\\
& u_{tt}\in L^2([0,T_*];L^2),\quad  t^{\frac{1}{2}}u\in L^\infty([0,T_*];D^4),\\
&t^{\frac{1}{2}}u_t\in L^\infty([0,T_*];D^2)\cap L^2([0,T_*] ; D^3),\ t^{\frac{1}{2}}u_{tt}\in L^\infty([0,T_*];L^2)\cap L^2([0,T_*];D^1).
\end{split}
\end{equation}
Moreover, if $1< \gamma \leq 3$, then $\rho(t,x)\in C^1([0,T_*]\times \mathbb{R}^2)$.
\end{theorem}

\begin{remark}\label{re3}
The smoothing effect of the velocity $u$ in positive time $t\in [\tau,T_*]$, $\forall \tau \in (0,T_*)$, tells us that the regular  solution obtained in Theorem \ref{th2} is indeed a classical one in $(0,T_*]\times \mathbb{R}^2$ (see \cite{guahu}).
\end{remark}

As a direct consequence of Theorem \ref{th2} and the standard theory of quasi-linear hyperbolic equations, we have
\begin{corollary}\label{co2} Let $1< \gamma \leq \frac{5}{3}$ or $\gamma=2,3$. If the initial data $( \rho_0, u_0)$ satisfy (\ref{th78}),
then there exists a time $T_*>0$ and a unique regular solution $(\rho, u)$ to the Cauchy problem  (\ref{eq:1.1}), satisfying (\ref{reg11}) and
\begin{equation}\label{regco}\begin{split}
& \rho\in C([0,T_*];H^3),\quad  \rho_t \in C([0,T_*];H^2),\quad \rho_{tt} \in C([0,T_*];L^2).
\end{split}
\end{equation}
\end{corollary}

We remark that the local existence time $T_*$ and all the estimates for the regularity of regular solutions in Theorem \ref{th2} depend only on  norms of $(\rho_0, u_0)$ stated in \eqref{th78}. The following theorem is about the continuous dependence of the solution on the initial data at least for a small time interval.

\begin{theorem}\label{co02} Let conditions in Theorem \ref{th2} hold. For each $i=1,2$, let $( \rho_i, u_i)$ be the local regular solutions to   Cauchy problem  (\ref{eq:1.1}) on $[0, T_{*}]$ with the initial data $( \rho_{0i}, u_{0i})$ satisfying (\ref{th78}). Let $K>0$ be a constant such that
$$
\left\|\rho^{\frac{\gamma-1}{2}}_{0i}\right\|_2+\left\|\frac{\nabla \rho_{0i}}{\rho_{0i}}\right\|_{L^6\cap D^1}+\|u_{0i}\|_2\leq K.
$$
Then for $0\leq t\leq T_*$,  there exists  a positive constant $C(T_*,K)$ such that
\begin{equation}\label{rezheng}\begin{split}
&\left\|\rho^{\frac{\gamma-1}{2}}_{1}-\rho^{\frac{\gamma-1}{2}}_{2}\right\|_2+\left\|\frac{\nabla \rho_{1}}{\rho_{1}}-\frac{\nabla \rho_{2}}{\rho_{2}}\right\|_{L^6\cap D^1}+\|u_{1}-u_{2}\|_2+\int_{0}^{T_*}\|\nabla u_1-\nabla u_2)\|^2_2\text{d}t\\
\leq &C\left(\left\|\rho^{\frac{\gamma-1}{2}}_{01}-\rho^{\frac{\gamma-1}{2}}_{02}\right\|_2+\left\|\frac{\nabla \rho_{01}}{\rho_{01}}-\frac{\nabla \rho_{01}}{\rho_{01}}\right\|_{L^6 \cap D^1}+\|u_{01}-u_{02}\|_2\right).
\end{split}
\end{equation}
\end{theorem}

Finally, we establish some blow-up criterion for classical solutions in terms of $\nabla \rho/\rho $ and $
D(u)=\frac{1}{2}\left(\nabla u+(\nabla u)^\top\right)$,
 which is analogous to the Beale-Kato-Majda criterion for the ideal incompressible flow \cite{pc}.

\begin{theorem}\label{th3s}
Let $(\rho(x,t), u(x,t))$ be a regular solution obtained in Theorem \ref{th2}.
If $\overline{T}< +\infty $ is the maximal existence time, then  both
\begin{equation}\label{cri}
\begin{split}
\displaystyle\lim_{T\mapsto \overline{T}} \left(\sup_{0\le t\le T}\left\|\frac{\nabla \rho}{\rho}(\cdot, t)\right\|_{L^6(\mathbb{R}^2)}+\int_0^T \|D( u)(\cdot, t)\|_ {L^\infty(\mathbb{R}^2)}\ dt\right)=+\infty,
\end{split}
\end{equation}
and
\begin{equation}\label{cri1}
\begin{split}
\lim \sup_{T\mapsto \overline{T}}\int_0^T \|D( u)(\cdot, t)\|_{L^\infty\cap D^{1,6}(\mathbb{R}^2)}\ dt=+\infty.
\end{split}
\end{equation}
\end{theorem}

As explained before, the main purposes of this article is to establish the theories applicable to most of physically relevant models in shallow water theory. Our framework for system \eqref{eq:1.1} is applicable with minor modifications to the viscous terms of the forms  $\text{div}(h D(U)), \ \text{div}(h\nabla U), \ h\Delta U$, and the one in \eqref{newvis}. More precisely, for viscous term in \eqref{newvis}, this is a special case of
system \eqref{eq:1.1} with $\alpha=1$, $\beta=2$ and $\gamma=2$.  When $\mathcal{V}=\text{div}(h D(U))$, this is also a special case of system \eqref{eq:1.1} with $\beta=0,\,\, \alpha=1/2,$ and $\gamma=2$.
Therefore, one simply replaces $\rho$ by $h$ in Theorems 1.1--1.3 to obtain the same results for these two classes of shallow water models, without further modifications. For viscous Saint-Venant model, $\mathcal{V}=\text{div}(h\nabla U)$, we let $\gamma=2$ in condition $(\textrm{A})$ and  $Q(U)=\nabla U$ in condition $(\textrm{D})$ of Definition \ref{d1}. Then we have the following theorem.

\begin{theorem}\label{th3}
If $\mathcal{V}=\text{div}(h\nabla U), $ and the initial data $( h_0, U_0)$ satisfy the regularity condition
\begin{equation}\label{th78q}
\begin{split}
&h_0> 0, \quad  (h_0, U_0)\in H^3, \quad \nabla h_0/h_0\in L^6\cap D^1\cap D^2,
\end{split}
\end{equation}
then there exists a time $T_*>0$ and a unique regular solution $(h, U)$ to Cauchy problem (\ref{eq:1.1ww}), satisfying
\begin{equation}\label{reg11q}\begin{split}
& h \in C([0,T_*];H^3),\ h_t \in C([0,T_*];H^2),\\
& \nabla h/h \in C([0,T_*]; L^6\cap D^1\cap D^2), \  (\nabla h/h)_t \in C([0,T_*]; H^1),\\
& U\in C([0,T_*]; H^3)\cap L^2([0,T_*] ; H^4), \ U_t \in C([0,T_*]; H^1)\cap L^2([0,T_*] ; D^2),\\
&U_{tt}\in L^2([0,T_*];L^2),\ t^{\frac{1}{2}}U\in L^\infty([0,T_*];D^4),\\
&
t^{\frac{1}{2}}U_t\in L^\infty([0,T_*];D^2)\cap L^2([0,T_*] ; D^3),\ t^{\frac{1}{2}}U_{tt}\in L^\infty([0,T_*];L^2)\cap L^2([0,T_*];D^1).
\end{split}
\end{equation}
Moreover, if  $\overline{T}< +\infty $ is the maximal existence time of the regular solution $(h, U)$, then
\begin{equation}\label{cri1w}
\begin{split}
\lim_{T\mapsto \overline{T}} \left(\sup_{0\le t\le T}\left\|\frac{\nabla h}{h}(\cdot, t)\right\|_{L^6(\mathbb{R}^2)}+\int_0^T\|\nabla U (\cdot, t)\|_{L^\infty(\mathbb{R}^2)}\right)=+\infty,
\end{split}
\end{equation}
and
\begin{equation}\label{cri2w}
\begin{split}
\lim_{T\mapsto \overline{T}}\left(\int_0^T \|\nabla U(\cdot, t)\|_{L^\infty\cap D^{1,6}(\mathbb{R}^2)}\ dt\right)=+\infty.
\end{split}
\end{equation}
\end{theorem}

Finally, when $\mathcal{V}=\ h\Delta U$, let $\gamma=2$ and $Q(U)=0$. Then we have the following result.

\begin{theorem}\label{zhengui}
If $\mathcal{V}=\ h\Delta U$, and the initial data $( h_0, U_0)$ satisfy the regularity condition
\begin{equation}\label{zhenth78q}
\begin{split}
&h_0\geq 0, \quad  (h_0, U_0)\in H^3,
\end{split}
\end{equation}
then there exists a time $T_*$ and a unique regular solution $(h, U)$ to  (\ref{eq:1.1ww}),  satisfying
\begin{equation}\label{zhen reg11q}\begin{split}
& h \in C([0,T_*];H^3),\quad h_t \in C([0,T_*];H^2),\\
& U\in C([0,T_*]; H^3)\cap L^2([0,T_*] ; H^4), \ U_t \in C([0,T_*]; H^1)\cap L^2([0,T_*] ; D^2),\\
& U_{tt}\in L^2([0,T_*];L^2),\quad  t^{\frac{1}{2}}U\in L^\infty([0,T_*];D^4),\\
&t^{\frac{1}{2}}U_t\in L^\infty([0,T_*];D^2)\cap L^2([0,T_*];D^3),\ t^{\frac{1}{2}}U_{tt}\in L^\infty([0,T_*];L^2)\cap L^2([0,T_*];D^1).
\end{split}
\end{equation}

\end{theorem}

\begin{remark}\label{shallow}
For this model, we note that the previous requirement on the regularity of $\nabla h_{0}/h_{0}$ is not needed due to the good structure of the viscous term $\mathcal {V}$. Current theorem allows initial data containing local vacuum.
\end{remark}

We now outline the organization of the rest of paper. In Section $2$, we list some important lemmas that will be used frequently in our proof. In Section $3$,  we first reformulate our problem  into a simpler form. Next we give the proof of the local existence of classical solutions to this reformulated problem. With the help of
a new variable $\psi=\frac{\nabla\rho}{\rho}$, and its symmetry as a gradient, we are able to achieve this local existence in four steps: 1) we construct approximate solutions for the linearized problem when initial density has positive lower bound; 2) we establish the a priori estimates independent of the lower bound of density for the linearized problem; 3) we then pass to the limit to recover the solution of this linearized problem allowing vacuum in the far field; 4) we prove the unique solvability of the reformulated problem through a standard iteration process. Section $4$ is devoted to proving the $H^2$ stability with respect to the initial data, i.e.
Theorem 1.2. The proof of Theorem 1.3 is given in Section 5.

Finally, we remark that our framework provided in this paper is applicable to the case of three dimensions, with some minor modifications. We will not pursue this in this article.

\section{Preliminaries}
In this section, we present some important lemmas  that appear frequently in our proof.
The first one is the following well-known Gagliardo-Nirenberg inequality, which can be found in \cite{oar}.
\begin{lemma}\cite{oar}\label{lem2as}
 Let $r\in (1,+\infty)$ and $ \ h\in W^{1,p}(\mathbb{R}^2) \cap L^r(\mathbb{R}^2) $. Then
$$
|h|_q\leq C|\nabla h|^\theta_p |h|^{1-\theta}_r
$$
where $\theta=\big(\frac{1}{r}-\frac{1}{q}\big)\big(\frac{1}{r}-\frac{1}{p}+\frac{1}{2}\big)^{-1}$. If $p< 2$,  then $q\in [r,\frac{2p}{2-p}]$ when $r<\frac{2p}{2-p} $; and $q\in [\frac{2p}{2-p},r]$ when $r\geq\frac{2p}{2-p} $. If $p=2$, then $q\in [r,+\infty)$. If $p> 2$, then  $q\in [r,+\infty]$.

\end{lemma}

Some common versions of this inequality can be written as
\begin{equation}\label{ine}\begin{split}
|f|_3\leq C|f|^{\frac{2}{3}}_{2}|\nabla f|^{\frac{1}{3}}_{2},\quad |f|_6\leq C|f|^{\frac{1}{3}}_{2}|\nabla f|^{\frac{2}{3}}_{2},\quad
|f|_\infty\leq C|f|^{\frac{1}{2}}_{6}|\nabla f|^{\frac{1}{2}}_{3},
\end{split}
\end{equation}
which will be used frequently in our following proof.

The second one  can be found in Majda \cite{amj}, and we omit its proof.

\begin{lemma}\cite{amj}\label{zhen1}
Let positive constants $r$, $a$ and $b$ satisfy the relation $\frac{1}{r}=\frac{1}{a}+\frac{1}{b}$ and  $1\leq a,\ b, \ r\leq +\infty$. $ \forall s\geq 1$, if $f, g \in W^{s,a} (\mathbb{R}^2)\cap  W^{s,b}(\mathbb{R}^2)$, then we have
\begin{equation}\begin{split}\label{ku11}
&|D^s(fg)-f D^s g|_r\leq C_s\big(|\nabla f|_a |D^{s-1}g|_b+|D^s f|_b|g|_a\big),
\end{split}
\end{equation}
\begin{equation}\begin{split}\label{ku22}
&|D^s(fg)-f D^s g|_r\leq C_s\big(|\nabla f|_a |D^{s-1}g|_b+|D^s f|_a|g|_b\big),
\end{split}
\end{equation}
where $C_s> 0$ is a constant depending on $s$ only.
\end{lemma}

The following lemma is important in the derivation of  uniqueness in Section $3$, which can be found in Remark 1 of \cite{bjr}.
\begin{lemma}\cite{bjr}\label{1}
If $f(t,x)\in L^2([0,T]; L^2)$, then there exists a sequence $s_k$ such that
$$
s_k\rightarrow 0, \quad \text{and}\quad s_k |f(s_k,x)|^2_2\rightarrow 0, \quad \text{as} \quad k\rightarrow+\infty.
$$
\end{lemma}

Due to  harmonic analysis,  we have the following regularity estimate result for Lam$\acute{ \text{e}}$ operator. For problem
\begin{equation}\label{ok}
-\alpha\triangle u-(\alpha+\beta)\nabla \text{div}u =Lu=F, \quad  u\rightarrow 0 \quad \text{as} \ |x|\rightarrow +\infty,
\end{equation}
we have
\begin{lemma}\cite{harmo}\label{zhenok}
If $u\in D^{1,q}(\mathbb{R}^2)$ with $1< q< +\infty$ is a weak solution to problem (\ref{ok}), then
$$
|u|_{D^{k+2,q}} \leq C |F|_{D^{k,q}},
$$
where $C$ depending only on $\alpha$,  $\beta$  and $q$.
\end{lemma}
The  proof can be obtained via the classical estimates from harmonic analysis, which can be found in  \cite{harmo} or \cite{zif}.

Finally, using Aubin-Lions Lemma, one has (c.f. \cite{jm}),
\begin{lemma}\cite{jm}\label{aubin} Let $X_0$, $X$ and $X_1$ be three Banach spaces satisfying $X_0\subset X\subset X_1$. Suppose that $X_0$ is compactly embedded in $X$ and that $X$ is continuously embedded in $X_1$.\\[1pt]

I) Let $G$ be bounded in $L^p(0,T;X_0)$ with $1\leq p < +\infty$, and $\frac{\partial G}{\partial t}$ be bounded in $L^1(0,T;X_1)$. Then $G$ is relatively compact in $L^p(0,T;X)$.\\[1pt]

II) Let $F$ be bounded in $L^\infty(0,T;X_0)$  and $\frac{\partial F}{\partial t}$ be bounded in $L^q(0,T;X_1)$ with $q>1$. Then $F$ is relatively compact in $C(0,T;X)$.
\end{lemma}

\section{Existence of Regular Solutions}

In this section, we aim at proving Theorem \ref{th2}. To this end, we first reformulated our main problem  (\ref{eq:1.1}) into a simpler form.

\subsection{Reformulation}
Let $\phi=\rho^{\frac{\gamma-1}{2}}$. System (\ref{eq:1.1})  can be written as
\begin{equation}
\begin{cases}
\label{eq:cccq}
\phi_t+u\cdot \nabla \phi+\frac{\gamma-1}{2}\phi\text{div} u=0,\\[12pt]
u_t+u\cdot\nabla u +\frac{2A\gamma}{\gamma-1}\phi\nabla \phi+Lu=\psi\cdot Q(u),
 \end{cases}
\end{equation}
where $\psi=\nabla \rho/\rho=\frac{2}{\gamma-1}\nabla \phi/\phi=(\psi^{(1)},\psi^{(2)})^\top$.
The initial data are given by
\begin{equation}\label{qwe}
(\phi,u)|_{t=0}=(\phi_0,u_0),\quad x\in \mathbb{R}^2.
\end{equation}

To prove Theorem \ref{th2}, our first step is to establish the following existence result for the reformulated problem (\ref{eq:cccq})-(\ref{qwe}).
\begin{theorem}\label{th1} If the initial data $( \phi_0,u_0)$ satisfy the following regularity conditions:
\begin{equation}\label{th78qq}
\begin{split}
&\phi_0> 0,\quad  (\phi_0, u_0)\in H^3, \quad  \psi_0=\frac{2}{\gamma-1}\nabla \phi_0/\phi_0\in  L^6\cap D^1\cap D^2,
\end{split}
\end{equation}
then there exist a time $T_*>0$ and a unique regular solution $( \phi,u)$ to Cauchy problem (\ref{eq:cccq})-(\ref{qwe}), satisfying
\begin{equation}\label{reg11qq}\begin{split}
& \phi \in C([0,T_*];H^3),\quad  \phi_t \in C([0,T_*];H^2),\\
&\psi \in C([0,T_*]; L^6\cap D^1\cap D^2),\ \psi_t \in C([0,T_*]; H^1),\\
& u\in C([0,T_*]; H^3)\cap L^2([0,T_*] ; H^4), \ u_t \in C([0,T_*]; H^1)\cap L^2([0,T_*] ; D^2),\\
& u_{tt}\in L^2([0,T_*];L^2),\quad  t^{\frac{1}{2}}u\in L^\infty([0,T_*];D^4),\\
&t^{\frac{1}{2}}u_t\in L^\infty([0,T_*];D^2)\cap L^2([0,T_*] ; D^3),\ t^{\frac{1}{2}}u_{tt}\in L^\infty([0,T_*];L^2)\cap L^2([0,T_*];D^1).
\end{split}
\end{equation}
\end{theorem}

We will prove this theorem in subsequent four subsections, and at the end of this section we will show that this theorem indeed implies Theorem \ref{th2}. For simplicity, in the following sections, we denote $\frac{A\gamma}{\gamma-1}=\theta$.

\subsection{Linearization}\label{linear2}

 In order to proceed with nonlinear problem, we first need to consider the following linearized problem
\begin{equation}\label{li4}
\begin{cases}
\phi_t+v\cdot \nabla \phi+\frac{\gamma-1}{2}\phi\text{div} v=0,\\[8pt]
u_t+v\cdot\nabla v +2\theta\phi\nabla \phi+Lu=\psi\cdot Q(v),\\[8pt]
(\phi,u)|_{t=0}=(\phi_0,u_0),\quad x\in \mathbb{R}^2,
 \end{cases}
\end{equation}
 where
\begin{equation} \label{eq:5.231}Lu=-\alpha\triangle u-(\alpha+\beta)\nabla \text{div}u,\ Q(v)=\alpha(\nabla v+(\nabla v)^\top)+\beta\text{div}v\mathbb{I}_{2}, \end{equation}
and $v=(v^{(1)},v^{(2)})\in \mathbb{R}^2$ is a known vector satisfying $v(t=0,x)=u_0(x)$ and
\begin{equation}\label{vg}
\begin{split}
& v\in C([0,T] ; H^3)\cap L^2([0,T] ; H^4), \ v_t \in C([0,T] ; H^1)\cap L^2([0,T] ; D^2),\\
& v_{tt}\in L^2([0,T];L^2),\ t^{\frac{1}{2}}v\in L^\infty([0,T];D^4),\\
& t^{\frac{1}{2}}v_t\in L^\infty([0,T];D^2)\cap L^2([0,T] ; D^3),\  t^{\frac{1}{2}}v_{tt}\in L^\infty([0,T];L^2)\cap L^2([0,T];D^1).
\end{split}
\end{equation}

We assume that
\begin{equation}\label{th78rr}
\begin{split}
&\phi_0> 0,\quad (\phi_0-\phi^\infty, u_0)\in H^3, \quad  \psi_0=\frac{2}{\gamma-1}\nabla \phi_0/\phi_0 \in  L^6\cap D^1\cap D^2,
\end{split}
\end{equation}
where $\phi^\infty\geq 0$ is a constant.

In the following two subsections, we first solve  this linearized problem when the initial density is away from vacuum, then we establish the uniform estimates with respect to the lower bound of the density which
enable us to pass to the limit of the case when $\inf \phi_{0}=0$.

\subsection{A priori estimate with uniformly positive density ($\inf \phi_{0}=\delta>0$).}

Now we want to get some local (in time) a priori estimate  which is independent of the lower bound of $\phi$ for the classical solution  $(\phi, u)$ to (\ref{li4}).
First we have  the following existence of classical solutions to (\ref{li4}) by the standard hyperbolic theory.

 \begin{lemma}\label{lem1}
 Assume that the initial data  $(\phi_{0}, u_{0})$ satisfy (\ref{th78rr}) and $\phi_0 \ge \delta$ for some constant $\delta>0$.
Then there exists a unique classical solution $(\phi,u)$ to (\ref{li4}) such that
\begin{equation}\label{reggh}\begin{split}
&\phi\geq \underline{\delta}, \ \phi-\phi^\infty \in C([0,T]; H^3), \ \phi _t \in C([0,T]; H^2), \\
&\psi \in C([0,T] ; L^6\cap D^1\cap D^2),\ \psi_t \in C([0,T]; H^1),\  \psi_{tt} \in  L^2([0,T]; L^2), \\
& u\in C([0,T]; H^3)\cap L^2([0,T] ; H^4), \ u_t \in C([0,T]; H^1)\cap L^2([0,T] ; D^2),\\
& u_{tt}\in L^2([0,T];L^2),\ t^{\frac{1}{2}}u\in L^\infty([0,T];D^4),\\
& t^{\frac{1}{2}}u_t\in L^\infty([0,T];D^2)\cap L^2([0,T] ; D^3),\ t^{\frac{1}{2}}u_{tt}\in L^\infty([0,T];L^2)\cap L^2([0,T];D^1),
\end{split}
\end{equation}
for any $T>0$, where $\underline{\delta}>0$ is a constant depending on $\delta$ and $T$.
\end{lemma}
\begin{proof}
The existence and regularity of a unique solution $\phi$ to the first equation of (\ref{li4}) can be obtained essentially according to Lemma 6 in  \cite{CK} via the standard theory of transport equation, and $\phi$ can be written as
\begin{equation}
\label{eq:bb1}
\phi(t,x)=\phi_0(W(0,t,x))\exp\Big(-\frac{\gamma-1}{2}\int_{0}^{t}\textrm{div}v(s,W(s,t,x))\text{d}s\Big),
\end{equation}
where  $W\in C^1([0,T]\times[0,T]\times \mathbb{R}^2)$ is the solution to the initial value problem
\begin{equation}
\label{eq:bb1}
\begin{cases}
\frac{d}{dt}W(t,s,x)=v(t,W(t,s,x)),\quad 0\leq t\leq T,\\[6pt]
W(s,s,x)=x, \quad \ \quad \quad 0\leq s\leq T,\ x\in \mathbb{R}^2.
\end{cases}
\end{equation}
So we easily know that there exists a positive constant $\underline{\delta}$ such that $\phi\geq \underline{\delta}$.

Now, it is easy to show that $\psi$ satisfies
$$
\psi_t+\nabla (v\cdot \psi)+\nabla \text{div} v=0.
$$
A direct calculation shows that $\partial_i \psi^{(j)}=\partial_j \psi^{(i)}$ in the sense of distribution,  then the above system can be rewritten as
\begin{equation}\label{ku}
\psi_t+\sum_{l=1}^2 A_l \partial_l\psi+B\psi+\nabla \text{div} v=0,
\end{equation}
where $A_l=(a^{(l)}_{ij})_{2\times 2}$ ($i,j,l=1,2$) are symmetric  with $a^{(l)}_{ij}=v^{(l)}$ when $i=j$ and $a^{(l)}_{ij}=0$ otherwise, $B=(\nabla v)^\top$. Therefore, system (\ref{ku}) is a positive symmetric system. Then the assertion of the regularity on $\psi$ follows.

Finally, with the regularity properties of $\phi$ and $\psi$, it is not difficult to solve $u$ from the linear parabolic equations
$$
u_t+Lu=-v\cdot\nabla v -2\theta\phi\nabla \phi+\psi\cdot Q(v),
$$
to complete the proof of this lemma. Here we omit the details.
\end{proof}

Now we are going to establish the uniform estimates on the solutions obtained in Lemma \ref{lem1}.
For this purpose, we fix $T>0$ and a positive constant $c_0$ large enough such that
\begin{equation}\label{houmian}\begin{split}
2+\phi^\infty+|\phi_0|_{\infty}+\|\phi_0-\phi^\infty\|_{3}+|\psi_0|_{L^6\cap D^1\cap D^2}+\|u_0\|_{3}\leq c_0.
\end{split}
\end{equation}
We now assume that there exist some time $T^*\in (0,T)$ and constants $c_i$ ($i=1,2,3,4$) such that
$$1< c_0\leq c_1 \leq c_2 \leq c_3 \leq c_4, $$
and
\begin{equation}\label{jizhu1}
\begin{split}
\sup_{0\leq t \leq T^*}\| v(t)\|^2_{1}+\int_{0}^{T^*} \Big( \|\nabla v(s)\|^2_{1}+|v_t(s)|^2_{2}\Big)\text{d}s \leq& c^2_1,\\
\sup_{0\leq t \leq T^*}\big(|v(t)|^2_{D^2}+|v_t(t)|^2_{2}\big)+\int_{0}^{T^*} \Big( |v(s)|^2_{D^3}+|v_t(s)|^2_{D^1}\Big)\text{d}s \leq& c^2_2,\\
\sup_{0\leq t \leq T^*}\big(|v(t)|^2_{D^3}+|v_t(t)|^2_{D^1}\big)+\int_{0}^{T^*} \Big( |v(s)|^2_{D^4}+|v_t(s)|^2_{D^2}+|v_{tt}(s)|^2_2\Big)\text{d}s \leq& c^2_3,\\
\text{ess}\sup_{0\leq t \leq T^*}\big(t|v_t(t)|^2_{D^2}+t|v(t)|^2_{D^4}+t|v_{tt}(t)|^2_{2}\big)+\int_{0}^{T^*} \Big(s|v_{tt}|^2_{D^1}+s|v_{t}|^2_{D^3}\Big)\text{d}s \leq& c^2_4.
\end{split}
\end{equation}
We shall determine $T^*$ and  $c_i$ ($i=1,2,3,4$) later, see \eqref{dingyi}, such that they
depend only on $c_0$ and the fixed constants $\alpha$, $\beta$, $\gamma$ and $T$.

Let $(\phi,u)$ be the unique classical solution to (\ref{li4}) on $[0,T] \times \mathbb{R}^2$. In the following we are going to establish a series of uniform local (in time) estimates listed as Lemmas  \ref{2}-\ref{5}. We start with the estimates for $\phi$.
Hereinafter, we use $C\geq 1$ to denote  a generic positive constant depending only on fixed constants $\alpha$, $\beta$, $\gamma$ and $T$.

\begin{lemma}\label{2} Let $(\phi,u)$ be the unique classical solution to (\ref{li4}) on $[0,T] \times \mathbb{R}^2$. Then
\begin{equation}\label{diyi}
\begin{split}
|\phi(t)|^2_\infty+\|\phi(t)-\phi^\infty\|^2_3\leq& Cc^2_0,\\
 |\phi_t(t)|_2\leq Cc_0c_1,\quad  |\phi_t(t)|_{D^1}\leq& Cc_0c_2,\\
 |\phi_t(t)|^2_{D^2}+|\phi_{tt}(t)|^2_2+\int_0^t \|\phi_{tt}(s)\|^2_1 \text{d}s\leq& Cc^6_3
\end{split}
\end{equation}
for $0\leq t \leq T_1=\min (T^{*}, (1+c_3)^{-2})$.
\end{lemma}

\begin{proof}
From stand energy estimates (see, for instance, \cite{CK}) and (\ref{ine}),  we easily have
\begin{equation}\label{gb}\begin{split}
\|\phi(t)-\phi^\infty\|_{3}\leq \Big(\|\phi_0-\phi^\infty\|_{3} +\phi^\infty\int_0^t \|\nabla v(s)\|_{3}\text{d}s\Big)\exp\Big(C\int_0^t \| v(s)\|_{4}\text{d}s\Big).
\end{split}
\end{equation}
Therefore, observing that
$$
\int_0^t \| v(s)\|_{4}\text{d}s\leq t^{\frac{1}{2}}\Big(\int_0^t \| v(s)\|^2_{4}\text{d}s\Big)^{\frac{1}{2}}\leq C(c_3t+c_3t^{\frac{1}{2}}),
$$
we get
$$\|\phi(t)-\phi^{\infty}\|_3 \leq C c_{0} \quad \text{for}\quad 0\leq t \leq T_1=\min (T^{*}, (1+c_3)^{-2}).$$

For $\phi_t$, since
$$\phi_t=-v\cdot \nabla \phi-\frac{\gamma-1}{2}\phi\text{div} v,$$
then for $0\leq t \leq T_1$, it holds that\begin{equation}\label{zhen6}
\begin{split}
|\phi_t(t)|_2\leq& C(|v\cdot \nabla \phi|_2+|\phi\text{div} v|_2)(t)\leq Cc_0c_1,\\
| \phi_t(t)|_{D^1}\leq& C\big(| v\cdot \nabla \phi|_{D^1}+| \phi\text{div} v|_{D^1}\big)(t)\leq Cc_0c_2,\\
| \phi_t(t)|_{D^2}\leq& C\big(| v\cdot \nabla \phi|_{D^2}+| \phi\text{div} v|_{D^2}\big)(t)\leq Cc_0c_3,
\end{split}
\end{equation}
where we used (\ref{ine}) and H\"older's inequality.

Using the equation
$$\phi_{tt}=-v_t\cdot \nabla \phi-v\cdot \nabla \phi_t-\frac{\gamma-1}{2}\phi_t\text{div} v-\frac{\gamma-1}{2}\phi\text{div} v_t$$
and the assumption (\ref{jizhu1}), we have,  for $0\leq t \leq T_1$,
\begin{equation}\label{zhen7}
\begin{split}
|\phi_{tt}(t)|_2
\leq&  C\big(|v_t\cdot \nabla \phi(t)|_2+|v\cdot \nabla \phi_t(t)|_2+|\phi_t\text{div} v(t)|_2+|\phi\text{div} v_t(t)|_2\big)\\
\leq& C\big(\|v_t(t)\|_1\|\phi(t)-\phi^\infty\|_3+\|  \phi_t(t)\|_1\|v(t)\|_3\big)\leq Cc^3_3.
\end{split}
\end{equation}

Similarly, for $0\leq t \leq T_1$, we also have
\begin{equation}\label{uu1}
\begin{split}
\int_0^t \|\phi_{tt}\|^2_1\text{d}s
\leq& C\int_0^t \big(\|v_t\cdot \nabla \phi\|^2_1+\|v\cdot \nabla \phi_t\|^2_1+\|\phi_t\text{div} v\|^2_1+\|\phi\text{div} v_t\|^2_1\big) \text{d}s\\
\leq& C\int_0^t \big(\|v_t\|^2_2\|\phi-\phi^\infty\|^2_3+\|v\|^2_3\|\phi_t\|^2_2\big) \text{d}s
\leq Cc^6_3.
\end{split}
\end{equation}
This completes the proof of the lemma.
\end{proof}

Now we establish the estimates for $\psi$ by the stand energy estimates for positive symmetric hyperbolic system.
\begin{lemma}\label{3} Let $(\phi,u)$ be the unique classical solution to (\ref{li4}) on $[0,T] \times \mathbb{R}^2$. Then
\begin{equation}\label{psi}
\begin{split}
&|\psi(t)|^2_\infty+\|\psi(t)\|^2_{L^6\cap D^1\cap D^2}\leq Cc^2_0,\quad |\psi_t|_2\leq Cc^2_2,\\
&|\psi_t(t)|^2_{D^1}+\int_0^t |\psi_{tt}(s)|^2_{2}\text{d}s\leq Cc^4_3,\quad \text{for}\quad  0\leq t \leq T_1.
\end{split}
\end{equation}
\end{lemma}
\begin{proof}
According to the proof of Lemma \ref{lem1}, we know that $\psi$  satisfies the  hyperbolic system (\ref{ku}).
First, multiplying (\ref{ku}) by $6|\psi|^4\psi$ and then integrating over $\mathbb{R}^2$, we easily deduce that
\begin{equation}\label{zhenzhen22}\begin{split}
\frac{d}{dt}| \psi|^6_6\leq C\Big(\sum_{l=1}^{2}|\partial_{l}A_l|_\infty+|B|_\infty\Big)|\psi|^6_6+C|\nabla^2 v|_6| \psi|^5_6.
\end{split}
\end{equation}
Noting that
\begin{equation}\label{zhenchen1}\begin{split}
|\nabla v|_\infty \leq C \|v\|_3,\
|\nabla^2 v|_6 \leq C|\nabla^2 v|^{\frac{1}{3}}_{2}|\nabla^3 v|^{\frac{2}{3}}_{2}\leq C\|\nabla^2 v\|_1,
\end{split}
\end{equation}
 (\ref{zhenzhen22}) implies that
\begin{equation}\label{zhenzhen33}\begin{split}
\frac{d}{dt}| \psi|_6\leq C\|v\|_3|\psi|_6+C\|\nabla^2 v\|_1.
\end{split}
\end{equation}

Second, let $\varsigma=(\varsigma_1,\varsigma_2)^\top$ ($1\leq |\varsigma|\leq 2$ and $\varsigma_i=0,1,2$). Taking derivative $\partial_{x}^{\varsigma} $ to (\ref{ku}), we have
\begin{equation}\label{hyp}\begin{split}
&(\partial_{x}^{\varsigma}  \psi)_t+\sum_{l=1}^2 A_l \partial_l\partial_{x}^{\varsigma}  \psi+B\partial_{x}^{\varsigma}  \psi+\partial_{x}^{\varsigma}  \nabla \text{div} v \\
=&\Big(-\partial_{x}^{\varsigma} (B\psi)+B\partial_{x}^{\varsigma}  \psi\Big)+\sum_{l=1}^2 \Big(-\partial_{x}^{\varsigma} (A_l \partial_l \psi)+A_l \partial_l\partial_{x}^{\varsigma}  \psi\Big)=\Theta_1+\Theta_2.
\end{split}
\end{equation}
Multiplying (\ref{hyp}) by $2\partial_{x}^{\varsigma} \psi$ and then integrating over $\mathbb{R}^2$, we have
\begin{equation}\label{zhenzhen}\begin{split}
\frac{d}{dt}|\partial_{x}^{\varsigma}  \psi|^2_2\leq \Big(\sum_{l=1}^{2}|\partial_{l}A_l|_\infty+|B|_\infty\Big)|\partial_{x}^{\varsigma}  \psi|^2_2+(|\Theta_1 |_2+|\Theta_2|_2+\|\nabla^2 v\|_2)|\partial_{x}^{\varsigma}  \psi|_2.
\end{split}
\end{equation}
For $|\varsigma|=1$, we apply (\ref{ku22}) with the choice $r=a=2$ and $b=+\infty$ to obtain,
\begin{equation}\label{zhen2}
|\Theta_1|_2=|\partial_{x}^{\varsigma} (B\psi)-B\partial_{x}^{\varsigma}  \psi|_2\leq C|\nabla^2 v|_2|\psi|_\infty;
\end{equation}
while the choice $r=b=2$, $a=\infty$  gives
\begin{equation}\label{zhen2l}
|\Theta_2|_2=|\partial_{x}^{\varsigma} (A_l \partial_l \psi)-A_l \partial_l\partial_{x}^{\varsigma}  \psi|_2\leq C|\nabla v|_\infty |\nabla\psi|_2.
\end{equation}

Similarly, for $|\varsigma|=2$, we have
\begin{equation}\label{zhen2}
|\Theta_1|_2=|\partial_{x}^{\varsigma} (B\psi)-B\partial_{x}^{\varsigma}  \psi|_2\leq C|\nabla^2 v|_\infty|\nabla\psi|_2+C|\nabla^3 v|_2 |\psi|_\infty,
\end{equation}
and
\begin{equation}\label{zhen2t}
|\Theta_2|_2=|\partial_{x}^{\varsigma} (A_l \partial_l \psi)-A_l \partial_l\partial_{x}^{\varsigma}  \psi|_2\leq C|\nabla v|_\infty|\nabla^2\psi|_2+C|\nabla^2 v|_\infty |\nabla\psi|_2.
\end{equation}
Using the fact that
\begin{equation}\label{ine1}\begin{split}
|\psi|_\infty\leq C|\psi|^{\frac{1}{2}}_{6}|\nabla \psi|^{\frac{1}{2}}_{3}\leq C|\psi|^{\frac{1}{2}}_{6} |\nabla \psi|^{\frac{1}{3}}_{2}|\nabla^2 \psi |^{\frac{1}{6}}_{2}\leq \|\psi\|_{L^6\cap D^1\cap D^2},
\end{split}
\end{equation}
formulas (\ref{zhenzhen22})-(\ref{ine1}) and Gagliardo-Nirenberg inequality lead to
\begin{equation*}
\frac{d}{dt}|\psi(t)|_{L^6\cap D^1\cap D^2}\leq C\| v \|_4\|\psi(t)\|_{L^6\cap D^1\cap D^2}+C\|\nabla ^2 v\|_2.
\end{equation*}
Then the Gronwall's inequality implies that
\begin{equation}\label{uu2}\begin{split}
|\psi(t)|_{L^6\cap D^1\cap D^2}\leq&  \Big(|\psi_0|_{L^6\cap D^1\cap D^2}+\int_0^t \|v\|_{4} \text{d}t\Big) \exp\Big(C\int_0^t \|v\|_{4} \text{d}t\Big)\leq Cc_0.
\end{split}
\end{equation}
For $0\leq t \leq T_1$, this estimate gives the part of  $\|\psi(t)\|_{L^6\cap D^1\cap D^2}$ in this lemma.

Noting that
\begin{equation}\label{ghtu}\psi_t=-\nabla (v \cdot \psi)-\nabla \text{div} v,
\end{equation}
for  $0\leq t \leq T_1$, it holds that
\begin{equation}\label{uu3}\begin{cases}
|\psi_t(t)|_2\leq C\big(|v|_{\infty}| \psi|_{D^1}+|\nabla v|_2|\psi|_{\infty}+|v|_{D^2}\big)(t)\leq Cc^2_2,\\[9pt]
|\nabla \psi_t(t)|_{2}\leq C \big(|\nabla v|_{\infty}|\psi|_{D^1}+|v|_{\infty}|\psi|_{D^2}+|\psi|_{\infty}|\nabla^2 v|_2+|v|_{D^3}\big)(t) \leq Cc^2_3.
\end{cases}
\end{equation}
Similarly, using
$$\psi_{tt}=-\nabla (v \cdot \psi)_t-\nabla \text{div} v_t,$$
for $0\leq t \leq T_1$,  we have
\begin{equation}\label{uu4}\begin{split}
\int_0^t |\psi_{tt}|^2_{2} \text{d}t
\leq& C\int_0^t \big(c^2_3|\nabla \psi|^2_3+c^2_3|\psi_t|^2_{2}+c^2_3|\nabla \psi_t|^2_2+c^2_3|\nabla v_t|^2_{2}+|v_t|^2_{D^2}\big) \text{d}t
\leq Cc^4_3.
\end{split}
\end{equation}
This concludes the proof of the lemma.
\end{proof}

Now we turn to the estimate of the velocity $u$.

 \begin{lemma}\label{4} Let $(\phi,u)$ be the unique classical solution to (\ref{li4}) on $[0,T] \times \mathbb{R}^2$. Then
\begin{equation}\label{uu}
\begin{split}
\| u(t)\|^2_{1}+\int_{0}^{T_2} \Big( \|\nabla u(s)\|^2_{1}+|u_t(s)|^2_{2}\Big)\text{d}s \leq& Cc^2_0,\\
|u(t)|^2_{D^2}+|u_t(t)|^2_{2}+\int_{0}^{T_2} \Big( |u(s)|^2_{D^3}+|u_t(s)|^2_{D^1}\Big)\text{d}s \leq& Cc^{\frac{10}{3}}_1c^{\frac{2}{3}}_2,
\end{split}
\end{equation}
for $0 \leq t \leq T_2=\min(T^*,(1+c_3)^{-6})$.
 \end{lemma}
\begin{proof} We divide the proof into three steps. \\
\underline{Step 1} (Estimate of $|u|_2$). Multiplying  $(\ref{li4})_2$ by $u$ and integrating over $\mathbb{R}^2$, we have
\begin{equation}\label{zhu1}
\begin{split}
&\frac{1}{2} \frac{d}{dt}|u|^2_2+\alpha|\nabla u|^2_2+(\alpha+\beta)|\text{div} u|^2_2\\
=&\int_{\mathbb{R}^2} \Big(-v\cdot \nabla v \cdot u-2\theta\phi \nabla \phi \cdot u+\psi \cdot Q(v)\cdot u \Big) \text{d}x\equiv: \sum_{i=1}^3 I_i.
\end{split}
\end{equation}
Due to Gagliardo-Nirenberg inequality, H\"older's inequality and Young's inequality, we have
\begin{equation}\label{zhu2}
\begin{split}
I_1=&-\int_{\mathbb{R}^2} v\cdot \nabla v \cdot u \text{d}x\leq C|v|_\infty|\nabla v|_2|u|_2\leq C\|v\|^2_2|\nabla v|^2_2+C| u|^2_2,\\
I_2=&-\int_{\mathbb{R}^2} 2\theta\phi \nabla \phi \cdot u \text{d}x
\leq C|\nabla\phi|_2| \phi|_\infty| u|_2\leq C| u|^2_2+C|\nabla\phi|^2_2|\phi|^2_\infty,\\
I_3=&\int_{\mathbb{R}^2} \psi \cdot Q(v)\cdot u  \text{d}x\leq C|\psi|_6|\nabla v|_3|u|_2\leq C|u|^2_2+C|\psi|^2_6\|\nabla v \|^2_1.
\end{split}
\end{equation}
Then we have
\begin{equation}\label{zhu3}
\begin{split}
&\frac{1}{2} \frac{d}{dt}|u|^2_2+\alpha|\nabla u|^2_2\leq C|u|^2_2+C\|v\|^2_2|\nabla v|^2_2+C|\nabla\phi|^2_2|\phi|^2_\infty+C|\psi|^2_6\|\nabla v \|^2_1.
\end{split}
\end{equation}
Integrating (\ref{zhu3}) over $(0,t)$, for $0 \leq t \leq T_1$, it gives
\begin{equation*}
\begin{split}
|u(t)|^2_2+\int_0^t\alpha|\nabla u(s)|^2_2\text{d}s\leq  C\int_0^t  |u(s)|^2_2 \text{d}s+C|u_0|^2_2+Cc^{4}_3t.
\end{split}
\end{equation*}
Then Gronwall's inequality implies
\begin{equation}\label{zhu5}
\begin{split}
|u(t)|^2_2+\int_0^t\alpha|\nabla u(s)|^2_2\text{d}s\leq  C\Big(|u_0|^2_2+c^{4}_3t\Big)\exp(Ct)\leq Cc^2_0
\end{split}
\end{equation}
for $0 \leq t \leq T_2=\min(T^*,(1+c_3)^{-6})$.

\underline{Step 2} (Estimate of $|\nabla u|_2$). Multiplying $ (\ref{li4})_2$ by $u_t$ and integrating over $\mathbb{R}^2$, we have
\begin{equation}\label{zhu6}
\begin{split}
&\frac{1}{2} \frac{d}{dt}\Big(\alpha|\nabla u|^2_2+(\alpha+\beta)|\text{div}u|^2_2\Big)+| u_t|^2_2\\
=&\int_{\mathbb{R}^2} \Big(-v\cdot \nabla v \cdot u_t-2\theta\phi \nabla \phi \cdot u_t+\psi \cdot Q(v)\cdot u_t \Big) \text{d}x\equiv: \sum_{i=4}^6 I_i.
\end{split}
\end{equation}
Using Gagliardo-Nirenberg inequality, H\"older's inequality and Young's inequality, we have
\begin{equation}\label{zhu10}
\begin{split}
I_4=&-\int_{\mathbb{R}^2} v\cdot \nabla v \cdot u_t \text{d}x\leq C|v|_\infty|\nabla v|_2|u_t|_2\leq C\|v\|^2_2|\nabla v|^2_2+\frac{1}{10}| u_t|^2_2,\\
I_5=&-\int_{\mathbb{R}^2} 2\theta\phi \nabla \phi \cdot u_t \text{d}x
\leq C|\phi|_\infty | \nabla \phi|_2| u_t|_2\leq \frac{1}{10}| u_t|^2_2+C|\nabla \phi|^2_2|\phi|^2_\infty,\\
I_6=&\int_{\mathbb{R}^2} \psi \cdot Q(v)\cdot u_t  \text{d}x\leq C|\psi|_6|\nabla v|_3|u_t|_2\leq \frac{1}{10}| u_t|^2_2+C|\psi|^2_6\|\nabla v\|^2_1.
\end{split}
\end{equation}
Then
\begin{equation}\label{zhu6q}
\begin{split}
&\frac{d}{dt}\Big(\alpha|\nabla u|^2_2+(\alpha+\beta)|\text{div}u|^2_2\Big)+| u_t|^2_2\\
\leq &C\|v\|^2_2|\nabla v|^2_2+C|\nabla \phi|^2_2|\phi|^2_\infty+C|\psi|^2_6\|\nabla v\|^2_1.
\end{split}
\end{equation}
Integrating (\ref{zhu6q}) over $(0,t)$, we get
\begin{equation}\label{zhu11}
\begin{split}
|\nabla u(t)|^2_2+\int_0^t| u_t(s)|^2_2\text{d}s\leq  C|\nabla u_0|^2_2+Cc^{4}_3t\leq Cc^2_0
\end{split}
\end{equation}
for $0 \leq t \leq T_2=\min(T^*,(1+c_3)^{-6})$.

 From the classical estimates for elliptic system in Lemma \ref{zhenok}, and
\begin{equation}\label{zhu77}
Lu=-u_t-v\cdot\nabla v -2\theta\phi\nabla \phi+\psi\cdot Q(v),
\end{equation}
 we easily have, for $0 \leq t \leq T_2$,
\begin{equation}\label{gai11}
\begin{split}
|u(t)|_{D^2}\leq& C \big(|u_t(t)|_2+|v\cdot\nabla v(t)|_2 +|\phi\nabla \phi(t)|_2+|\psi\cdot Q(v)(t)|_2\big)\\
\leq &C\Big(|u_t(t)|_2+c^{\frac{5}{3}}_1c^{\frac{1}{3}}_2+c^2_0+c_0|\nabla v|_2\Big),
\end{split}
\end{equation}
where we have used the fact that
$$
|v\cdot\nabla v(t)|_2\leq C|v|_6|\nabla v|_3\leq C|v|^{\frac{1}{3}}_2|\nabla v|^{\frac{4}{3}}_2|\nabla^2 v|^{\frac{1}{3}}_2.
$$
Then (\ref{gai11}) implies that
\begin{equation}\label{zhu77gai}
\int_0^{T_2}|u(t)|^2_{D^2} \text{d}t\leq C\int_0^{T_2} \Big(|u_t|^2_2+c^{\frac{10}{3}}_1c^{\frac{2}{3}}_2+c^4_0+c^2_0|\nabla v|^2_2\Big)(t)\text{d}t\leq Cc^2_0.
\end{equation}

\underline{Step 3} (Estimate of $|u|_{D^2}$). First we differentiate $(\ref{li4})_3$ with respect to $t$:
\begin{equation}\label{zhu7}
\begin{split}
u_{tt}+Lu_t=-(v\cdot\nabla v)_t -2\theta(\phi\nabla \phi)_t+(\psi\cdot Q(v))_t.
\end{split}
\end{equation}
Multiplying (\ref{zhu7}) by $u_t$ and integrating over $\mathbb{R}^2$, we have
\begin{equation}\label{zhu8}
\begin{split}
&\frac{1}{2} \frac{d}{dt}|u_t|^2_2+\alpha|\nabla u_t|^2_2+(\alpha+\beta)|\text{div} u_t|^2_2\\
=&\int_{\mathbb{R}^2} \Big(-(v\cdot \nabla v)_t \cdot u_t-2\theta(\phi \nabla \phi)_t \cdot u_t+(\psi \cdot Q(v))_t\cdot u_t \Big) \text{d}x\equiv: \sum_{i=7}^9 I_i,
\end{split}
\end{equation}
where the right-hand side terms can be estimated as follows,
\begin{equation}\label{zhu9}
\begin{split}
I_7=&-\int_{\mathbb{R}^2} (v\cdot \nabla v)_t \cdot u_t \text{d}x\leq C|v|_{\infty}|\nabla v_t|_{2}|u_t|_2+C|v_t|_2 |\nabla v |_\infty |u_t|_2\\
\leq& C\|v\|^2_{3}\| v_t\|^2_{1}+C|u_t|^2_2,\\
I_8=&-\int_{\mathbb{R}^2} 2\theta(\phi \nabla \phi)_t \cdot u_t \text{d}x=\theta\int_{\mathbb{R}^2}   (\phi^2)_t \text{div}u_t \text{d}x\\
\leq& C|\phi_t|_2| \phi|_\infty|\nabla u_t|_2\leq \frac{\alpha}{10}|\nabla u_t|^2_2+C|\phi_t|^2_2|\phi|^2_\infty,\\
I_9=&\int_{\mathbb{R}^2} (\psi \cdot Q(v))_t\cdot u_t\text{d}x\leq C|\psi|_\infty|\nabla v_t|_2|u_t|_2+C|\psi_t|_2 |\nabla v|_\infty |u_t|_2\\
\leq&C|u_t|^2_{2}+C|\psi|^2_{L^6\cap D^1\cap D^2}|\nabla v_t|^2_2+C|\psi_t|^2_2\|\nabla v\|^2_2.
\end{split}
\end{equation}
Then we have
\begin{equation}\label{zhu8qq}
\begin{split}
&\frac{d}{dt}|u_t|^2_2+|\nabla u_t|^2_2\\
\leq& C|u_t|^2_2+C\|v\|^2_{3}\|v_t\|^2_{1}+C|\phi_t|^2_2|\phi|^2_\infty+C|\psi|^2_{L^6\cap D^1\cap D^2}|\nabla v_t|^2_2+C|\psi_t|^2_2\|\nabla v\|^2_2.
\end{split}
\end{equation}

Integrating (\ref{zhu8qq}) over $(\tau,t)$ $(\tau \in( 0,t))$, we have
\begin{equation}\label{zhu13}
\begin{split}
&|u_t(t)|^2_2+\int_\tau^t |\nabla u_t(s)|^2_2 \text{d}s\\
\leq & |u_t(\tau)|^2_2+C\int_{\tau}^t  \Big(|u_t|^2_2+\|v\|^2_{3}\|v_t\|^2_{1}+|\phi_t|^2_2\|\phi\|^2_2\Big)(s)\text{d}s\\
&+C\int_{\tau}^t  (|\psi|^2_{L^6\cap D^1\cap D^2}|\nabla v_t|^2_2+|\psi_t|^2_2\|\nabla v\|^2_2)(s)\text{d}s\\
\leq & |u_t(\tau)|^2_2+C\int_0^t  |u_t(s)|^2_2 \text{d}s+Cc^6_3t,\quad  \text{for} \  0 \leq t \leq T_2.
\end{split}
\end{equation}

From the momentum equations $ (\ref{li4})_2$ we have
\begin{equation}\label{zhu15}
\begin{split}
|u_t(\tau)|_2\leq C\big( |v|_\infty |\nabla v|_2+|\phi|_\infty|\nabla \phi|_2+|u|_{D^2}+|\psi|_6|\nabla v|_2\big)(\tau).
\end{split}
\end{equation}
Then it is clear from the assumption (\ref{vg}) and Lemma \ref{lem1} that
\begin{equation}\label{zhu15vb}
\begin{split}
\lim \sup_{\tau\rightarrow 0}|u_t(\tau)|_2\leq C\big( |v_0|_\infty |\nabla v_0|_2+|\phi_0|_\infty|\nabla \phi_0|_2+|u_0|_{D^2}+|\psi_0|_6|\nabla v_0|_2\big)\leq Cc^2_0.
\end{split}
\end{equation}

Letting $\tau\rightarrow 0$ in (\ref{zhu13}), it follows from Gronwall's inequality that
\begin{equation}\label{zhu14}
\begin{split}
&|u_t(t)|^2_2+\int_0^t\Big(\alpha|\nabla u_t|^2_2+(\alpha+\beta)|\text{div} u_t|^2_2\Big)(s)\text{d}s
\leq  c^4_0\exp\big(Ct\big)\leq Cc^4_0
\end{split}
\end{equation}
for $0 \leq t \leq T_2$.

So from (\ref{gai11}), we easily have
 for $0 \leq t \leq T_2$,
\begin{equation*}
\begin{split}
|u(t)|_{D^2}\leq &C\Big(|u_t(t)|_2+c^{\frac{5}{3}}_1c^{\frac{1}{3}}_2+c^2_0+c_0|\nabla v|_2\Big)\leq Cc^{\frac{5}{3}}_1c^{\frac{1}{3}}_2.
\end{split}
\end{equation*}
From the classical estimates for elliptic system in Lemma \ref{zhenok}, we have
\begin{equation}\label{jiabiao}
\begin{split}
|u(s)|_{D^3}\leq& C \big(|u_t|_{D^1}+|v\cdot\nabla v|_{D^1} +|\phi\nabla \phi|_{D^1}+|\psi\cdot Q(v)|_{D^1}\big)\\
\leq &C\big(|u_t|_{D^1}+|v|_6|\nabla^2 v|_3+|\nabla v|_6|\nabla v|_3+|\phi|_\infty|\nabla^2 \phi|_2\big)\\
&+C\big(|\nabla \phi|_6|\nabla \phi|_3+|\psi|_\infty|\nabla^2 v|_2+|\nabla\psi|_6|\nabla v|_3\big)\\
\leq &C\Big(|u_t|_{D^1}+c_1c_2+c_1c^{\frac{2}{3}}_2c^{\frac{1}{3}}_3\Big),
\end{split}
\end{equation}
where we have used the face that
$$
|v|_6|\nabla^2 v|_3\leq C|v|^{\frac{1}{3}}_2|\nabla v|^{\frac{2}{3}}_2|\nabla^2 v|^{\frac{2}{3}}_2|\nabla^3 v|^{\frac{1}{3}}_2,\quad |\nabla v|_6|\nabla v|_3\leq C|\nabla v|^{\frac{1}{3}}_2|\nabla^2 v|^{\frac{2}{3}}_2|\nabla v|^{\frac{2}{3}}_2|\nabla^2 v|^{\frac{1}{3}}_2.
$$
Then (\ref{jiabiao}) implies that
$$
\int_0^{T_2} |u(s)|^2_{D^3}\text{d}t\leq C\int_0^{T_2}\Big(|u_t|^2_{D^1}+c^4_3\Big) \text{d}t\leq Cc^4_0.
$$
\end{proof}

Now we will give some estimates for the higher order terms of the velocity $u$ in the following  Lemma.
\begin{lemma}\label{5} Let $(\phi,u)$ be the unique classical solution to (\ref{li4}) on $[0,T] \times \mathbb{R}^2$. Then,
\begin{equation}\label{lipan1}
\begin{split}
&|u_t(t)|^2_{D^1}+|u(t)|^2_{D^3}+\int_{0}^{t}\Big(|u_t(s)|^2_{D^{2}}+|u_{tt}(s)|^2_{2}+|u(s)|^2_{D^{4}}\Big)\text{d}s\leq Cc^{\frac{10}{3}}_2c^{\frac{2}{3}}_3, \\
&t|u_t(t)|^2_{ D^2}+t|u_{tt}(t)|^2_2+t|u(t)|^2_{D^4}+\int_{0}^{t}(s|u_{tt}|^2_{D^1}+s|u_{t}|^2_{D^3})\text{d}s\leq Cc^{4}_3
\end{split}
\end{equation}
for $0\leq t \leq T_3=\min (T^*, (1+c_3)^{-8})$.
 \end{lemma}
\begin{proof}We divide the proof into two steps.

\underline{Step 1} (Estimate of $|u|_{D^3}$).
Multiplying  (\ref{zhu7}) by $u_{tt}$ and integrating over $\mathbb{R}^2$,  we have
\begin{equation}\label{zhu19}
\begin{split}
& \frac{1}{2}\frac{d}{dt}\Big(\alpha|\nabla u_t|^2_2+(\alpha+\beta)|\text{div}u_t|^2_2\Big)+| u_{tt}|^2_2\\
=&\int_{\mathbb{R}^2} \Big(-(v\cdot \nabla v)_t\cdot u_{tt}-2\theta(\phi \nabla \phi)_t \cdot u_{tt}+(\psi \cdot Q(v))_t\cdot u_{tt} \Big) \text{d}x\equiv: \sum_{i=10}^{12} I_i.
\end{split}
\end{equation}

Applying Gagliardo-Nirenberg inequality, H\"older's inequality and Young's inequality, we get
\begin{equation}\label{zhu20}
\begin{split}
I_{10}=&-\int_{\mathbb{R}^2} (v\cdot \nabla v)_t \cdot u_{tt} \text{d}x\leq C|v_t|_{2}|\nabla v|_{\infty}|u_{tt}|_2+C|v|_\infty |\nabla v_t|_2 |u_{tt}|_2\\
\leq& C|v_t|^2_{2}\|\nabla v\|^2_{2}+C\|v\|^2_2 |\nabla v_t|^2_2 +\frac{1}{10}|u_{tt}|^2_2,\\
I_{11}=&-\int_{\mathbb{R}^2} 2\theta(\phi \nabla \phi)_t \cdot u_{tt} \text{d}x\leq C|\phi_t|_2 |\nabla \phi|_\infty |u_{tt}|_2+C|\phi|_\infty |\nabla \phi_t|_2 |u_{tt}|_2\\
\leq& C|\phi_t|^2_2 \|\nabla \phi\|^2_2 +C|\phi|^2_\infty |\nabla \phi_t|^2_2 +\frac{1}{10}|u_{tt}|^2_2,\\
\end{split}
\end{equation}
\begin{equation}\label{jiajiagaigai}
\begin{split}
I_{12}=&\int_{\mathbb{R}^2} (\psi \cdot Q(v))_t\cdot u_{tt}\text{d}x\leq C|\psi_t|_2|\nabla v|_\infty|u_{tt}|_2+C|\psi|_\infty |\nabla v_t|_2 |u_{tt}|_2\\
\leq& C|\psi_t|^2_2\|\nabla v\|^2_2 +C|\psi|^2_{L^6\cap D^1\cap D^2} |\nabla v_t|^2_2  +\frac{1}{10}|u_{tt}|^2_2.
\end{split}
\end{equation}
Then
\begin{equation}\label{zhu19qq}
\begin{split}
& \frac{d}{dt}\Big(\alpha|\nabla u_t|^2_2+(\alpha+\beta)|\text{div}u_t|^2_2\Big)+| u_{tt}|^2_2\\
\leq &C|v_t|^2_{2}\|\nabla v\|^2_{2}+C\|v\|^2_2 |\nabla v_t|^2_2+C|\phi_t|^2_2 \|\nabla \phi\|^2_2 +C|\phi|^2_\infty |\nabla \phi_t|^2_2\\
&+C|\psi_t|^2_2\|\nabla v\|^2_2 +C|\psi|^2_{L^6\cap D^1\cap D^2} |\nabla v_t|^2_2.
\end{split}
\end{equation}
Integrating (\ref{zhu19qq}) over $(\tau,t)$,  we have
\begin{equation}\label{zhu22g}
\begin{split}
|\nabla u_t(t)|^2_2+\int_{\tau}^t| u_{tt}(s)|^2_2\text{d}s\leq  C|\nabla u_t(\tau)|^2_2+Cc^{6}_3t,
\end{split}
\end{equation}
for $0\leq t \leq T_3=\min (T^*, (1+c_3)^{-8})$.

On the other hand, from the momentum equations $ (\ref{li4})_2$ we have
\begin{equation}\label{zhu15wsx}
\begin{split}
|\nabla u_t(\tau)|_2\leq C\big(1+ \|v\|^2_3+\|\phi\|^2_3+|\psi|_{L^6\cap D^1\cap D^2}\| v\|_2\big)(\tau).
\end{split}
\end{equation}
Then from the assumption (\ref{vg}) and Lemma \ref{lem1}, one has
\begin{equation}\label{zhu15vbwsx}
\begin{split}
\lim \sup_{\tau\rightarrow 0}|\nabla u_t(\tau)|_2\leq C\big(1+ \|v_0\|^2_3+\|\phi_0\|^2_3+|\psi_0|_{L^6\cap D^1\cap D^2}\| v_0\|_2\big)\leq Cc^2_0.
\end{split}
\end{equation}
Letting $\tau \rightarrow 0$ in (\ref{zhu22g}), it reads that
\begin{equation}\label{zhu22wsx}
\begin{split}
|\nabla u_t(t)|^2_2+\int_{0}^t| u_{tt}(s)|^2_2\text{d}s\leq  Cc^4_0,\quad \text{for}\ 0\leq t\leq T_3.
\end{split}
\end{equation}

For the higher order terms, from (\ref{jiabiao}) and (\ref{zhu22wsx}), it is easy to show that
\begin{equation*}
\begin{split}
|u(t)|_{D^3}\leq& C \Big(|u_t(t)|_{D^1}+c_1c_2+c_1c^{\frac{2}{3}}_2c^{\frac{1}{3}}_3\Big)\leq Cc^{\frac{5}{3}}_2c^{\frac{1}{3}}_3.
\end{split}
\end{equation*}
From (\ref{zhu77}) we have
\begin{equation}\label{zhu77qq}
Lu_t=-u_{tt}-(v\cdot\nabla v)_t -2\theta(\phi\nabla \phi)_t+(\psi\cdot Q(v))_t.
\end{equation}
We apply Lemma \ref{zhenok} to \eqref{zhu77qq} to show that, for $0 \leq t \leq T_3$,
\begin{equation}\label{gaijia1}
\begin{split}
|u_t(t)|_{D^2}\leq& C\big(|u_{tt}|_2+|(v\cdot\nabla v)_t|_2 +|(\phi\nabla \phi)_t|_2+|(\psi\cdot Q(v))_t|_2\big)(t)\\
\leq &C\Big(|u_{tt}|_2+c_2|\nabla v_t|_2+c^{\frac{4}{3}}_2|\nabla v_t|^{\frac{2}{3}}_2+c^3_2+c^2_2|\nabla^3 v|_2\Big)(t),\\
|u(t)|_{D^4}\leq &C \big(|u_t|_{D^2}+|v\cdot\nabla v|_{D^2} +|\phi\nabla \phi|_{D^2}+|\psi\cdot Q(v)|_{D^2}\big)(t)\\
\leq &C \big(|u_t|_{D^2}+c^2_0+c_2\|\nabla v\|_2\big)(t),
\end{split}
\end{equation}
which quickly implies that
$$
\int_0^{T_3}\big( |u_t(t)|^2_{D^2}+ |u(t)|^2_{D^4}\big)\text{d}t\leq Cc^4_0,\quad \text{for} \quad 0 \leq t \leq T_3.
$$
\underline{Step 2} (Estimate of $t^{\frac{1}{2}}|u|_{D^4}$).
Now we differentiate (\ref{zhu7}) with respect to $t$:
\begin{equation}\label{zhu25}
\begin{split}
&u_{ttt}+Lu_{tt}
=-2v_t\cdot\nabla v_t -v_{tt}\cdot \nabla v-v\cdot \nabla v_{tt}\\
&-\theta(\nabla \phi^2)_{tt}+2\psi_t\cdot (Q(v))_t+\psi_{tt}\cdot Q(v)+ \psi\cdot(Q(v))_{tt}.
\end{split}
\end{equation}
Multiplying  (\ref{zhu25}) by $u_{tt}$ and integrating over $\mathbb{R}^2$,  we have
\begin{equation}\label{zhu27}
\begin{split}
&\frac{1}{2} \frac{d}{dt}|u_{tt}|^2_2+\alpha|\nabla u_{tt}|^2_2+(\alpha+\beta)|\text{div} u_{tt}|^2_2\\
=&\int_{\mathbb{R}^2} \Big(-2v_t\cdot\nabla v_t -v_{tt}\cdot \nabla v-v\cdot \nabla v_{tt}-\theta(\nabla \phi^2)_{tt}+2\psi_t\cdot (Q(v))_t\Big)\cdot u_{tt}  \text{d}x\\
&+\int_{\mathbb{R}^2} \Big(\psi_{tt}\cdot Q(v)+ \psi\cdot(Q(v))_{tt} \Big)\cdot u_{tt}  \text{d}x\equiv: \sum_{i=13}^{19} I_i.
\end{split}
\end{equation}
Similarly, we can estimate the right-hand side term by term as follows.
\begin{equation}\label{zhu26}
\begin{split}
I_{13}=&-\int_{\mathbb{R}^2} 2v_t\cdot\nabla v_t \cdot u_{tt} \text{d}x\leq C| u_{tt}|_2|\nabla v_t|_3|v_t|_6
\leq  C| u_{tt}|^2_2+C\| v_t\|^2_1\|\nabla v_t\|^2_1,\\
I_{14}=&-\int_{\mathbb{R}^2} v_{tt}\cdot \nabla v  \cdot u_{tt} \text{d}x
\leq C|\nabla v|_\infty | v_{tt}|_2| u_{tt}|_2
\leq C| u_{tt}|^2_2+C\|\nabla v\|^2_2 | v_{tt}|^2_2.
\end{split}
\end{equation}
For the term $I_{15}$, via the integration by parts, we have
\begin{equation}\label{zhu26gaigai}
\begin{split}
I_{15}=&-\int_{\mathbb{R}^2} v\cdot \nabla v_{tt}\cdot u_{tt}\text{d}x
\leq C|v|_\infty| v_{tt}|_2 |\nabla u_{tt}|_2+C|\nabla v|_\infty| v_{tt}|_2 | u_{tt}|_2\\
\leq&  \frac{\alpha}{20}| \nabla u_{tt}|^2_2+C| u_{tt}|^2_2+C\|v\|^2_3|v_{tt}|^2_2.
\end{split}
\end{equation}
Similarly, for the terms $I_{16}$-$I_{18}$, we have
\begin{equation}\label{zhu26qq}
\begin{split}
I_{16}=&\int_{\mathbb{R}^2} \theta (\phi^2)_{tt} \cdot \text{div}u_{tt}\text{d}x
\leq  C|\phi_{tt}|_2 |\phi|_\infty |\nabla u_{tt}|_2+ |\phi_t|_6 |\phi_t|_3 |\nabla u_{tt}|_2\\
\leq &\frac{\alpha}{10}|\nabla u_{tt}|^2_2+C|\phi_{tt}|^2_2|\phi|^2_\infty+C\|\phi_t\|^4_1,\\
I_{17}=&\int_{\mathbb{R}^2}  2\psi_t\cdot (Q(v))_t\cdot u_{tt}\text{d}x\leq C|\psi_t|_6|\nabla v_t|_3| u_{tt}|_2
\leq C| u_{tt}|^2_2+C\|\psi_t\|^2_1\|\nabla v_t\|^2_1,\\
I_{18}=&\int_{\mathbb{R}^2}  \psi_{tt}\cdot Q(v)\cdot u_{tt}\text{d}x\leq C| \psi_{tt}|_2|\nabla v|_\infty| u_{tt}|_2
\leq C| u_{tt}|^2_2+C| \psi_{tt}|^2_2\|\nabla v\|^2_2.
\end{split}
\end{equation}
For the last term, via the integration by parts, we have
\begin{equation}\label{zhu26gai}
\begin{split}
I_{19}=&\int_{\mathbb{R}^2}  \psi\cdot(Q(v))_{tt}\cdot u_{tt}\text{d}x
\leq C|\psi|_\infty|v_{tt}|_2| \nabla u_{tt}|_2+C|\nabla \psi|_3|v_{tt}|_2|  u_{tt}|_6\\
\leq& \frac{\alpha}{20}| \nabla u_{tt}|^2_2+C|\psi|^2_{L^6\cap D^1\cap D^2}| v_{tt}|^2_2+C| u_{tt}|^2_2.
\end{split}
\end{equation}
These estimates, together with Lemmas \ref{2}-\ref{4} and Step 1 lead to
\begin{equation}\label{zhu28}
\begin{split}
& \frac{d}{dt}|u_{tt}|^2_2+|\nabla u_{tt}|^2_2\\
\leq& C| u_{tt}|^2_2+C\| v_t\|^2_1\|\nabla v_t\|^2_1+C|\phi_{tt}|^2_2 |\phi|^2_\infty+C\|\phi_t\|^4_1\\
&+C(\|v\|^2_3+\|\psi\|^2_{L^6\cap D^1\cap D^2})| v_{tt}|^2_2+C\| \psi_t\|^2_1\|\nabla v_t\|^2_1+C| \psi_{tt}|^2_2\|\nabla v\|^2_2\\
\leq & C| u_{tt}|^2_2+Cc^2_3| v_{tt}|^2_2+Cc^4_3|v_t|_{D^2}+Cc^2_3|\psi_{tt}|^2_2+Cc^{8}_3
\end{split}
\end{equation}
for $0\leq t \leq T_3$.

Multiplying both sides of (\ref{zhu28}) by $t$  and  integrating over $(\tau,t)$, we have
\begin{equation}\label{zhu29}
\begin{split}
&t|u_{tt}(t)|^2_2+\int_\tau^ts|\nabla u_{tt}(s)|^2_2\text{d}s\\
\leq& \tau|u_{tt}(\tau)|^2_2+\int_\tau^t |u_{tt}(s)|^2_2\text{d}s+ C\int_\tau^t s\Big(| u_{tt}(s)|^2_2+c^2_3| v_{tt}(s)|^2_2\Big)\text{d}s\\
&+ C\int_\tau^t s\Big( c^4_3|v_t|_{D^2}+c^2_3|\psi_{tt}|^2_2\Big)\text{d}s+Cc^{8}_3t\\
\leq & \tau|u_{tt}(\tau)|^2_2+Cc^4_0+Cc^{8}_3t.
\end{split}
\end{equation}
We know from \eqref{lipan1} and Lemma \ref{1} that,
there exists a sequence $s_k$ such that
$$
s_k\rightarrow 0, \quad \text{and}\quad s_k |u_{tt}(s_k,x)|^2_2\rightarrow 0, \quad \text{as} \quad k\rightarrow+\infty.
$$
Taking $\tau=s_k$ and letting $ k\rightarrow+\infty$ in (\ref{zhu29}), we arrive at
\begin{equation}\label{zhu30}
\begin{split}
&t|u_{tt}(t)|^2_2+\int_\tau^ts\Big(\alpha|\nabla u_{tt}(s)|^2_2+(\alpha+\beta)|\text{div} u_{tt}(s)|^2_2\Big)\text{d}s\leq Cc^{4}_0.
\end{split}
\end{equation}

Now combining Lemmas \ref{2}-\ref{4}, (\ref{lipan1}) and (\ref{zhu77qq})-(\ref{gaijia1}),  from Lemma \ref{zhenok}, we obtain
\begin{equation*}
\begin{split}
&t^{\frac{1}{2}}|u_t(t)|_{D^2}+
t^{\frac{1}{2}}|u(t)|_{D^4}\leq Cc^{2}_3,\\
&t^{\frac{1}{2}}|u_t(t)|_{D^3}\leq Ct^{\frac{1}{2}}\big(c^4_3+|u_{tt}(t)|_{D^1}+c_3|v_t(t)|_{D^2} \big)
\end{split}
\end{equation*}
for $0\leq t\leq T_3$, which,  according to (\ref{zhu30}),  proves $(\ref{lipan1})_2$.
\end{proof}

Then from Lemmas \ref{2}-\ref{5}, for $0 \leq t \leq T_3$, we have
\begin{equation}\label{jkk11}
\begin{split}
\| u(t)\|^2_{1}+\int_{0}^{t} \Big(\ |\nabla u(s)\|^2_{1}+|u_t(s)|^2_{2}\Big)\text{d}s \leq& Cc^2_0,\\
|u(t)|^2_{D^2}+|u_t(t)|^2_{2}+\int_{0}^{t} \Big( |u(s)|^2_{D^3}+|u_t(s)|^2_{D^1}\Big)\text{d}s \leq& Cc^{\frac{10}{3}}_1c^{\frac{2}{3}}_2,\\
|u(t)|^2_{D^3}+|u_t(t)|^2_{D^1}+\int_{0}^{t} \Big( |u(s)|^2_{D^4}+|u_t(s)|^2_{D^2}+|u_{tt}(s)|^2_2\Big)\text{d}s \leq& Cc^{\frac{10}{3}}_2c^{\frac{2}{3}}_3,\\
t|u_t(t)|^2_{D^2}+t|u(t)|^2_{D^4}+t|u_{tt}(t)|^2_{2}+\int_{0}^{t} (s|u_{tt}|^2_{D^1}+s|u_{t}|^2_{D^3})\text{d}s \leq& Cc^4_3,\\
|\phi(t)|^2_\infty+\|\phi(t)-\phi^\infty\|^2_3+\|\phi_t(t)\|^2_2+|\phi_{tt}(t)|^2_2+\int_0^{t} \|\phi_{tt}\|^2_1 \text{d}s\leq& Cc^6_3,\\
|\psi(t)|^2_\infty+|\psi(t)|^2_{L^6\cap D^1\cap D^2}+\|\psi_t(t)\|^2_{1}+\int_0^{t} |\psi_{tt}|^2_{2}\text{d}s\leq& Cc^4_3.
\end{split}
\end{equation}
Therefore, if we define the constants $c_i$ ($i=1,2,3,4,5$) and $T^*$ by
\begin{equation}\label{dingyi}
\begin{split}
&c_1=C^{\frac{1}{2}}c_0, \quad  c_2=C^{\frac{3}{4}}c^{\frac{5}{2}}_1=C^2c^{\frac{5}{2}}_0,\quad c_3=C^{\frac{3}{4}}c^{\frac{5}{2}}_2=C^{\frac{23}{4}}c^{\frac{25}{4}}_0,\\
& c_4= C^{\frac{1}{2}}c^3_3=C^{\frac{71}{4}}c^{\frac{75}{4}}_0, \quad \text{and} \quad T^*=\min (T, (1+c_3)^{-8}),
\end{split}
\end{equation}
then we deduce that

\begin{equation}\label{jkk}
\begin{split}
\sup_{0\leq t \leq T^*}\| u(t)\|^2_{1}+\int_{0}^{T^*} \Big( \|\nabla u(s)\|^2_{1}+|u_t(s)|^2_{2}\Big)\text{d}s \leq& c^2_1,\\
\sup_{0\leq t \leq T^*}(|u(t)|^2_{D^2}+|u_t(t)|^2_{2})+\int_{0}^{T^*} \Big( |u(s)|^2_{D^3}+|u_t(s)|^2_{D^1}\Big)\text{d}s \leq& c^2_2,\\
\sup_{0\leq t \leq T^*}(|u(t)|^2_{D^3}+|u_t(t)|^2_{D^1})+\int_{0}^{T^*} \Big( |u(s)|^2_{D^4}+|u_t(s)|^2_{D^2}+|u_{tt}(s)|^2_2\Big)\text{d}s \leq& c^2_3,\\
ess\sup_{0\leq t \leq T^*}(t|u_t(t)|^2_{D^2}+t|u(t)|^2_{D^4}+t|u_{tt}(t)|^2_{2})+\int_{0}^{T^*} \Big(s|u_{tt}|^2_{D^1}+s|u_{t}|^2_{D^3}\Big)\text{d}s \leq& c^2_4,\\
\sup_{0\leq t \leq T^*}(|\phi(t)|^2_\infty+\|\phi(t)-\phi^\infty\|^2_3+\|\phi_t(t)\|^2_2+|\phi_{tt}(t)|^2_2)+\int_0^{T^*}  \|\phi_{tt}\|^2_1 \text{d}s\leq& c^2_4,\\
\sup_{0\leq t \leq T^*}(|\psi(t)|^2_\infty+|\psi(t)|^2_{L^6\cap D^1\cap D^2}+\|\psi_t(t)\|^2_{1})+\int_0^{T^*} |\psi_{tt}|^2_{2}\text{d}s\leq& c^2_4.
\end{split}
\end{equation}
In other word, given fixed $c_0$, $T$, there are positive constants $T^*$, $c_i$ $(i=1, 2, 3, 4)$, depending solely on $c_0$ and $T$, such that if \eqref{jizhu1} holds for $v$, then \eqref{jkk} holds for classical solution of
\eqref{li4} on $[0, T^*]\times \mathbb{R}^2$. We remark that this property is important for the iteration to be successful in the proof of Theorem 3.1 in section 3.5 below.

\subsection{Unique solvability with vacuum in far field}\ \\

Based on the local (in time) estimates in (\ref{jkk}), we have the following existence result under the assumption that $\phi_0> 0$.
\begin{lemma}\label{lem1q}
 Assume that the initial data $(\phi_{0}, u_{0})$ satisfy (\ref{th78qq}).
Then there exists a unique classical solution $(\phi,u)$ to (\ref{li4}) such that
\begin{equation}\label{regghq}\begin{split}
&\phi\geq 0, \ \phi \in C([0,T^*]; H^3), \ \phi _t \in C([0,T^*]; H^2),\ \psi \in C([0,T^*]; L^6\cap D^1\cap D^2), \\
&\partial_i \psi^{(j)}=\partial_j \psi^{(i)} \ (i,j=1,2),\ \psi_t \in C([0,T^*]; H^1),\  \psi_{tt} \in  L^2([0,T^*]; L^2), \\
& u\in C([0,T^*]; H^3)\cap L^2([0,T^*] ; H^4), \ u_t \in C([0,T^*]; H^1)\cap L^2([0,T^*] ; D^2),\\
& u_{tt}\in L^2([0,T^*];L^2),\ t^{\frac{1}{2}}u\in L^\infty([0,T^*];D^4),\\
& t^{\frac{1}{2}}u_t\in L^\infty([0,T^*];D^2)\cap L^2([0,T^*];D^3),\ t^{\frac{1}{2}}u_{tt}\in L^\infty([0,T^*];L^2)\cap L^2([0,T^*];D^1).
\end{split}
\end{equation}
Moreover, the solution $(\phi,u)$ also satisfies the estimates in (\ref{jkk}).
\end{lemma}
\begin{proof} For $\delta\in (0,1)$, we define
$$\phi_{\delta0}=\phi_0+\delta,\,\, \psi_{\delta0}=\frac{2}{\gamma-1}\nabla \phi_0/(\phi_0+\delta).$$

From \eqref{th78rr} we know that $\psi_0\in L^6\cap D^1\cap D^2(\mathbb{R}^2)$, then Gagliardo-Nirenberg inequality implies that, there exists a finite and positive constant $\widetilde{C}$ such that
$$|\nabla \phi_0(x)|\leq \widetilde{C}\phi_0(x), \quad  \text{for} \quad x\in \mathbb{R}^2.$$
 So if $\phi_0(x)=0$, we immediately have $\nabla \phi_0(x)=0$, which means that
\begin{equation*}\begin{split}
\psi_{\delta0}=0,\ \text{if} \  \phi_0(x)=0; \quad \psi_{\delta0}-\psi_0= -\frac{\delta}{\delta+\phi_0} \psi_0,\ \text{if} \  \phi_0(x)>0.
\end{split}
\end{equation*}

From the assumption (\ref{houmian}), there exists a $\delta_{1}>0$ such that if $0<\delta<\delta_{1}$, then
\begin{equation*}\begin{split}
1+|\phi_{\delta0}|_\infty+\|\phi_{\delta0}-\delta\|_{3}+|\psi_{\delta0}|_{L^6\cap D^1\cap D^2}+\|u_0\|_{3}\leq Cc^2_0= \overline{c}_0.
\end{split}
\end{equation*}
Therefore, taking $(\phi_{\delta0}, u_0)$ as the initial data, problem \eqref{li4} admits a unique classical solution $(\phi^\delta,u^\delta)$ satisfying the local estimates in (\ref{jkk}). We note that the estimates in (\ref{jkk}) are independent of $\delta$, then there exists a subsequence (still denoted by $(\phi^\delta,\psi^\delta,u^\delta)$) converges to a limit $(\phi,\psi,u)$ in weak or weak* sense:
\begin{equation}\label{ruojixian}
\begin{split}
(\phi^\delta-\delta,u^\delta)\rightharpoonup  (\phi,u) \quad &\text{weak* \ in } \ L^\infty([0,T^*];H^3(\mathbb{R}^2)),\\
\psi^\delta  \rightharpoonup  \psi \quad &\text{weak* \ in } \ L^\infty([0,T^*];L^6\cap D^1\cap D^2(\mathbb{R}^2)),\\
\phi^\delta_t\rightharpoonup  \phi_t \quad &\text{weak* \ in } \ L^\infty([0,T^*];H^2(\mathbb{R}^2)),\\
 (\psi^\delta_t, u^\delta_t)\rightharpoonup  (\psi_t,u_t) \quad &\text{weak* \ in } \ L^\infty([0,T^*];H^1(\mathbb{R}^2)).
\end{split}
\end{equation}

In addition, for any $R> 0$, due to the Aubin-Lions Lemma  (see \cite{jm}) (i.e., Lemma \ref{aubin}), there exists a subsequence (still denoted by $(\phi^\delta,\psi^\delta,u^\delta)$) satisfying
\begin{equation}\label{ert}\begin{split}
(\phi^\delta,\psi^\delta,u^\delta)\rightarrow (\phi,\psi,u) \quad \text{in } \ C([0,T^*];H^1(B_R)),
\end{split}
\end{equation}
where $B_R$ is a ball centered at origin with radius $R$.
It is clear that $(\phi,\psi,u) $ also satisfies the local estimates in (\ref{jkk}).
So it is easy to show that $(\phi,u) $ is a weak solution of problem (\ref{li4}) satisfying the  regularity:
\begin{equation*}\begin{split}
&\phi> 0, \ \phi \in L^\infty([0,T^*]; H^3), \ \phi _t \in L^\infty([0,T^*]; H^2), \\
& \psi \in L^\infty([0,T^*]; L^6\cap D^1\cap D^2), \ \psi_t \in L^\infty([0,T^*]; H^1),\  \psi_{tt} \in  L^2([0,T^*]; L^2), \\
& u\in L^\infty([0,T^*]; H^3)\cap L^2([0,T^*] ; H^4), \ u_t \in L^\infty([0,T^*]; H^1)\cap L^2([0,T^*] ; D^2),\\
& u_{tt}\in L^2([0,T^*];L^2),\ t^{\frac{1}{2}}u\in L^\infty([0,T^*];D^4),\\
& t^{\frac{1}{2}}u_t\in L^\infty([0,T^*];D^2)\cap L^2([0,T^*];D^3),\ t^{\frac{1}{2}}u_{tt}\in L^\infty([0,T^*];L^2)\cap L^2([0,T^*];D^1).
\end{split}
\end{equation*}
We remark that, in this step, even though vacuum appears in the far field, $\psi$ satisfies $\partial_i \psi^{(j)}=\partial_j \psi^{(i)}$  $(i,j=1,2)$ and  solves the following positive symmetric hyperbolic system in the sense of distribution:
\begin{equation}\label{zhenzheng}
\psi_t+\nabla (v\cdot \psi)+\nabla \text{div} v=0.
\end{equation}

The uniqueness and time continuity for $(\phi,u)$ can be obtained by standard procedure and we omit the details here.
\end{proof}

\subsection{Proof of Theorem \ref{th1}}\ \\

Our proof  is based on the classical iteration scheme and the existence results for the linearized problem  in  Sections 3.2-3.4. Like in Section 3.3, we define constants $c_{0}$, $c_{1}$, $c_{2}$, $c_3$, $c_4$ and assume that
\begin{equation*}\begin{split}
&2+|\phi_0|_{\infty}+\|(\phi_0,u_0)\|_{3}+\|\psi_0\|_{L^6\cap D^1\cap D^2} \leq c_0.
\end{split}
\end{equation*}
Denote by $u^0\in C([0,T^*];H^3)\cap  L^2([0,T^*];H^4) $ the solution of the following Cauchy problem of heat equation:
\begin{equation*}
\begin{cases}
h_t-\triangle h=0 \quad \text{in} \quad (0,+\infty)\times \mathbb{R}^2,\\[8pt]
h(0)=u_0 \quad \text{in} \quad \mathbb{R}^2.
 \end{cases}
\end{equation*}
Using the regularity of $u^0=h$, we can always choose a time $T_{**}\in (0,T^*]$ such that
\begin{equation}\label{jizhu}
\begin{split}
\sup_{0\leq t \leq T_{**}}\| u^0(t)\|^2_{1}+\int_{0}^{T_{**}} \Big( \|\nabla u^0(s)\|^2_{1}+|u^0_t(s)|^2_{2}\Big)\text{d}s \leq& c^2_1,\\
\sup_{0\leq t \leq T_{**}}(|u^0(t)|^2_{D^2}+|u^0_t(t)|^2_{2})+\int_{0}^{T_{**}} \Big( |u^0(s)|^2_{D^3}+|u^0_t(s)|^2_{D^1}\Big)\text{d}s \leq& c^2_2,\\
\sup_{0\leq t \leq T_{**}}(|u^0(t)|^2_{D^3}+|u^0_t(t)|^2_{D^1})+\int_{0}^{T_{**}} \Big( |u^0(s)|^2_{D^4}+|u^0_t(s)|^2_{D^2}+|u^0_{tt}(s)|^2_2\Big)\text{d}s \leq& c^2_3,\\
\text{ess}\sup_{0\leq t \leq T_{**}}\big(t|u^0_t(t)|^2_{D^2}+t|u^0(t)|^2_{D^4}\big)+\int_{0}^{T_{**}}\Big( s|u^0_{tt}|^2_{D^1}+s|u^0_{t}|^2_{D^3}\Big)\text{d}s \leq& c^2_4.
\end{split}
\end{equation}

We first prove the existence of regular solutions.
Let $v=u^0$,  we can get a classical solution $(\phi^1, u^1)$ of problem (\ref{li4})  as well as function $\psi^{1}$. Inductively, we construct approximate sequences $(\phi^{k+1}, \psi^{k+1}, u^{k+1})$ as follows: given $(\phi^{k},\psi^k, u^{k})$ for $k\geq 1$, define $(\phi^{k+1}, \psi^{k+1}, u^{k+1})$  by solving the following problem:

\begin{equation}\label{li6}
\begin{cases}
\phi^{k+1}_t+u^{k}\cdot \nabla \phi^{k+1}+\frac{\gamma-1}{2}\phi^{k+1}\text{div} u^{k}=0,\\[8pt]
\psi^{k+1}_t+\sum_{l=1}^2 A_l(u^k) \partial_l\psi^{k+1}+B(u^k)\psi^{k+1}+\nabla \text{div}u^k=0,\\[8pt]
u^{k+1}_t+u^{k}\cdot\nabla u^{k} +2\theta\phi^{k+1}\nabla \phi^{k+1}=-Lu^{k+1}+\psi^{k+1}\cdot Q(u^{k}),\\[8pt]
(\phi^{k+1}, \psi^{k+1}, u^{k+1})|_{t=0}=(\phi_0,\psi_0,u_0).
 \end{cases}
\end{equation}
This problem was obtained from \eqref{li4} by replacing $v$ with $ u^{k}$. Then we know that $(\phi^{k}, \psi^{k}, u^{k})$ $(k=1,2,...)$  satisfy the estimates in (\ref{jkk}).

Next we are going to prove that the whole sequence $(\phi^k,\psi^k, u^k)$ converges strongly to a limit $(\phi,\psi, u)$ which satisfies the regularity (\ref{reg11qq}).
Let
\begin{equation*}
\overline{\phi}^{k+1}=\phi^{k+1}-\phi^k,\ \overline{\psi}^{k+1}=\psi^{k+1}-\psi^k,\ \overline{u}^{k+1}=u^{k+1}-u^k.
\end{equation*}
From (\ref{li6}) we have
 \begin{equation}
\label{eq:1.2w}
\begin{cases}
\ \  \displaystyle \overline{\phi}^{k+1}_t+u^k\cdot \nabla\overline{\phi}^{k+1} +\overline{u}^k\cdot\nabla\phi ^{k}+\frac{\gamma-1}{2}(\overline{\phi}^{k+1} \text{div}u^k +\phi ^{k}\text{div}\overline{u}^k)=0,\\[8pt]
\ \ \displaystyle \overline{\psi}^{k+1}_t+\sum_{l=1}^2 A_l(u^k) \partial_l\overline{\psi}^{k+1}+B(u^k)\overline{\psi}^{k+1}+\nabla \text{div}\overline{u}^k=\Upsilon^k_1+\Upsilon^k_2,\\[8pt]
\ \ \overline{u}^{k+1}_t+ u^k\cdot\nabla \overline{u}^{k}+ \overline{u}^{k} \cdot \nabla u^{k-1}+\theta\nabla ((\phi^{k+1})^2-(\phi^k)^2) \\[8pt]
\qquad=-L\overline{u}^{k+1}+\psi^{k+1}\cdot Q(\overline{u}^k)+\overline{\psi}^{k+1}\cdot Q(u^{k-1}),
\end{cases}
\end{equation}
where $\Upsilon^k_1$  and $\Upsilon^k_2$ are defined by
\begin{equation*}
\Upsilon^k_1=-\sum_{l=1}^2(A_l(u^k) \partial_l\psi^{k}-A_l(u^{k-1}) \partial_l\psi^{k}),\quad \Upsilon^k_2=-(B(u^k) \psi^{k}-B(u^{k-1}) \psi^{k}).
\end{equation*}

We first estimate $\|\overline{\phi}^{k+1}\|_1$. Multiplying $(\ref{eq:1.2w})_1$ by $2\overline{\phi}^{k+1}$ and integrating over $\mathbb{R}^2$, we have
\begin{equation*}
\begin{split}
\frac{d}{dt}|\overline{\phi}^{k+1}|^2_2=& -2\int_{\mathbb{R}^2}\Big(u^k\cdot \nabla\overline{\phi}^{k+1} +\overline{u}^k\cdot\nabla\phi ^{k}+\frac{\gamma-1}{2}(\overline{\phi}^{k+1} \text{div}u^k +\phi ^{k}\text{div}\overline{u}^k)\Big)\overline{\phi}^{k+1} \text{d}x,\\
\leq& C|\nabla u^k|_\infty|\overline{\phi}^{k+1}|^2_2+C |\overline{\phi}^{k+1}|_2|\overline{u}^k|_2|\nabla \phi^k|_\infty+C|\overline{\phi}^{k+1}|_2|\nabla\overline{u}^k|_2| \phi^k|_\infty,
\end{split}
\end{equation*}
which means that ($0<\eta \leq \frac{1}{10}$ is a constant)
\begin{equation}\label{go64}\begin{split}
\displaystyle
&\frac{d}{dt}|\overline{\phi}^{k+1}(t)|^2_2\leq A^k_\eta(t)|\overline{\phi}^{k+1}(t)|^2_2+\eta \|\overline{u}^k(t)\|^2_1,\\[4pt]
\displaystyle
&A^k_\eta(t)=C\left(\|\nabla u^k\|_{2}+\frac{1}{\eta}(1+\|\phi^{k}\|^2_{3})\right),\ \text{and} \ \int_0^t A^k_\eta(s)\text{d}s\leq C +C_{\eta}t
\end{split}
\end{equation}
for $t\in[0,T_{**}]$, where $C_{\eta}$ is a positive constant  depending on $\eta$ and  constant $C$.

Now taking derivative $\partial_{x}^{\zeta}$ ($|\zeta|=1$) to $(\ref{eq:1.2w})_{1}$, multiplying by $2 \partial_{x}^\zeta\overline{\phi}^{k+1}$ and  integrating over $\mathbb{R}^2$,  we have
\begin{equation*}
\begin{split}
\frac{d}{dt}|\partial_{x}^\zeta\overline{\phi}^{k+1}|^2_{2}=&-2\int_{\mathbb{R}^2}\partial_{x}^\zeta\big(u^k\cdot \nabla\overline{\phi}^{k+1} +\overline{u}^k\cdot\nabla\phi ^{k})\big)\partial_{x}^\zeta\overline{\phi}^{k+1} \text{d}x\\
&-2\int_{\mathbb{R}^2}\partial_{x}^\zeta\big(\frac{\gamma-1}{2}(\overline{\phi}^{k+1} \text{div}u^k +\phi ^{k}\text{div}\overline{u}^k)\big)\partial_{x}^\zeta\overline{\phi}^{k+1} \text{d}x\\
\leq& C|\nabla u^k|_\infty |\nabla \overline{\phi}^{k+1}|^2_{2}+C|\nabla \phi^k|_{\infty}| \nabla \overline{u}^k|_2|\nabla \overline{\phi}^{k+1}|_2\\
&+C|\nabla \overline{\phi}^{k+1}|_{2}|  \overline{u}^k|_6|\nabla^2 \phi^{k}|_3+C|\nabla^2 u^k|_\infty | \overline{\phi}^{k+1}|_{2} |\nabla \overline{\phi}^{k+1}|_{2}\\
&+C|\nabla \overline{u}^k|_2 |\nabla \overline{\phi}^{k+1}|_{2}|\nabla \phi^{k}|_\infty+C| \phi^k|_{\infty}  |\nabla \text{div} \overline{u}^k|_2   |\nabla \overline{\phi}^{k+1}|_{2},
\end{split}
\end{equation*}
which means that
\begin{equation}\label{go64qqqq}
\begin{split}
\displaystyle
&\frac{d}{dt}|\nabla \overline{\phi}^{k+1}(t)|^2_2\leq B^k_\eta(t)\|\overline{\phi}^{k+1}(t)\|^2_1+\eta \|\overline{u}^k(t)\|^2_2,\\[4pt]
&B^k_\eta(t)=C\Big(\| u^k\|_{4}+\frac{1}{\eta}\|\phi^{k}\|^2_{3} \Big),\ \text{and} \ \int_0^t B^k_\eta(s)\text{d}s\leq C+C_{\eta}t
\end{split}
\end{equation}
for $t\in[0,T_{**}]$. Combining (\ref{go64})-(\ref{go64qqqq}), it is easy to show that,  for $t\in[0,T_{**}]$,
\begin{equation}\label{go64qqqqqq}
\begin{cases}
\displaystyle
\frac{d}{dt}\| \overline{\phi}^{k+1}(t)\|^2_1\leq \Phi^k_\eta(t)\| \overline{\phi}^{k+1}(t)\|^2_1+\eta \|\overline{u}^k(t)\|^2_2,\\[8pt]
\displaystyle \Phi^{k_{\eta}}(t)=C(A^{k}_{\eta}(t)+B^{k}_{\eta}(t)), \ \text{and} \
 \int_0^t \Phi^k_\eta(s)\text{d}s\leq C+C_{\eta}t.
\end{cases}
\end{equation}

Now we estimate $|\overline{\psi}^{k+1}|_2$. Multiplying $(\ref{eq:1.2w})_2$ by $2\overline{\psi}^{k+1}$ and integrating over $\mathbb{R}^2$, we have
\begin{equation}\label{go64aa}
\begin{split}
\frac{d}{dt}|\overline{\psi}^{k+1}|^2_2\leq& \Big(\sum_{l=1}^{2}|\partial_{l}A_l(u^k)|_\infty+|B(u^k)|_\infty\Big)|\overline{\psi}^{k+1}|^2_2\\
&+(|\Upsilon^k_1 |_2+|\Upsilon^k_2|_2+|\nabla^2 \overline{u}^k|_2)|\overline{\psi}^{k+1}|_2.
\end{split}
\end{equation}
From H\"older's inequality, it is easy to deduce that
\begin{equation}\label{go64aa1}
\begin{split}|\Upsilon^k_1 |_2\leq C|\nabla \psi^{k}|_3 |\overline{u}^k|_6, \quad  |\Upsilon^k_2 |_2\leq C| \psi^{k}|_\infty | \nabla \overline{u}^k|_2.
\end{split}
\end{equation}
From (\ref{go64aa})-(\ref{go64aa1}), for $t\in[0,T_{**}]$, we have
\begin{equation}\label{go64aa2}\begin{cases}
\displaystyle
\frac{d}{dt}|\overline{\psi}^{k+1}|^2_2\leq \Psi^k_\eta(t)|\overline{\psi}^{k+1}|^2_2+\eta \|\overline{u}^k\|^2_2,\\[10pt]
\displaystyle
\Psi^k_\eta(t)=C\Big(\|\nabla u^k\|_{2}+\frac{1}{\eta}|\psi^{k}|^2_{\infty} +\frac{1}{\eta}\ |\nabla\psi^{k}\ |^2_{3}+\frac{1}{\eta}\Big),\ \text{and} \ \int_0^t \Psi^k_\eta(s)\text{d}s\leq C+C_{\eta}t.
\end{cases}
\end{equation}

For $\|\overline{u}^{k+1}\|_1$, multiplying $(\ref{eq:1.2w})_3$ by $2\overline{u}^{k+1}$ and integrating over $\mathbb{R}^2$, we have
\begin{equation*}\begin{split}
&\frac{d}{dt}|\overline{u}^{k+1}|^2_2+2\alpha|\nabla\overline{u}^{k+1} |^2_2+2(\alpha+\beta)|\text{div}\overline{u}^{k+1} |^2_2\\
=& -2\int_{R^3}\Big( u^k \cdot \nabla \overline{u}^{k}+\overline{u}^{k} \cdot \nabla u^{k-1}\Big) \cdot \overline{u}^{k+1}\text{d}x\\
&-2\int_{R^3}\Big(\theta\nabla \big((\phi^{k+1})^2-(\phi^k)^2\big)-\psi^{k+1}\cdot Q(\overline{u}^k)+\overline{\psi}^{k+1}\cdot Q(u^{k-1}) \Big) \cdot \overline{u}^{k+1}\text{d}x\\
\leq &C|u^k|_\infty |\nabla\overline{u}^{k}|_2 |\overline{u}^{k+1}|_2+C|\overline{u}^{k}|_6 |\overline{u}^{k+1}|_2 |\nabla u^{k-1}|_3 +C|\psi^{k+1}|_\infty |\nabla\overline{u}^{k}|_2 |\overline{u}^{k+1}|_2 \\
&+C \Big(|\phi^{k+1}|_\infty+|\phi^{k}|_\infty\Big)|\nabla \overline{u}^{k+1}|_2 |\overline{\phi}^{k+1}|_2+C|\overline{\psi}^{k+1}|_2 |\nabla u^{k-1}|_\infty |\overline{u}^{k+1}|_2,
\end{split}
\end{equation*}
which implies that
\begin{equation}\label{gogo1}\begin{split}
&\frac{d}{dt}|\overline{u}^{k+1}|^2_2+\alpha |\nabla\overline{u}^{k+1} |^2_2\\
\leq& E^k_\eta(t)|\overline{u}^{k+1}|^2_2+E^k_2(t)|\overline{\phi}^{k+1}|^2_{2}
+E^k_3(t)|\overline{\psi}^{k+1}|^2_{2}+\eta\|\overline{u}^{k}\|^2_1,
\end{split}
\end{equation}
where
\begin{equation*}
\begin{split}
E^k_\eta(t)=&C\left(1+\frac{1}{\eta}|u^{k}|^2_{\infty}+\frac{1}{\eta}| \nabla u^{k-1}|^2_{3}+\frac{1}{\eta}|\psi^{k+1}|^2_{\infty}\right), \\ E^k_2(t)=&C(|\phi^{k+1}|_\infty+|\phi^{k}|_\infty)^2,\
E^k_3(t)=C|\nabla u^{k-1}|^2_\infty ,
\end{split}
\end{equation*}
and
$$\int_0^t \big(E^k_\eta(s)+E^k_2(s)+E^k_3(s)\big)\text{d}s\leq C+C_\eta t.$$

Now taking $\partial_{x}^{\zeta}$ to $(\ref{eq:1.2w})_3$ ($|\zeta|=1$), multiplying by $\partial_{x}^\zeta\overline{u}^{k+1}$ and  integrating over $\mathbb{R}^2$,  we have
\begin{equation*}
\begin{split}
&\frac{1}{2}\frac{d}{dt}|\partial_{x}^\zeta\overline{u}^{k+1}|^2_{2}+\alpha|\nabla \partial_{x}^\zeta\overline{u}^{k+1} |^2_2+(\alpha+\beta)| \partial_{x}^\zeta \text{div} \overline{u}^{k+1} |^2_2\\
=& \int_{\mathbb{R}^2} \Big(-\partial_{x}^\zeta( u^k\cdot\nabla \overline{u}^{k})-\partial_{x}^\zeta(\overline{u}^{k} \cdot \nabla u^{k-1})-\theta \partial_{x}^\zeta(\nabla ((\phi^{k+1})^2-(\phi^k)^2) )\Big)\cdot \partial_{x}^\zeta\overline{u}^{k+1} \text{d}x \\
&+\int_{\mathbb{R}^2} \Big(\partial_{x}^\zeta\big(\psi^{k+1}\cdot Q(\overline{u}^k)\big)+\partial_{x}^\zeta\big(\overline{\psi}^{k+1}\cdot Q(u^{k-1})\big)\Big)\cdot \partial_{x}^\zeta\overline{u}^{k+1} \text{d}x=\sum_{i=1}^5 J_i,
\end{split}
\end{equation*}
where
\begin{equation}
\begin{split}\label{litong3}
J_1=&\int_{\mathbb{R}^2} -\partial_{x}^\zeta( u^k\cdot\nabla \overline{u}^{k})\cdot \partial_{x}^\zeta\overline{u}^{k+1} \text{d}x\\
\leq &C|\nabla u^k|_\infty |\nabla\overline{u}^{k}|_2 |\nabla \overline{u}^{k+1}|_2+C| u^k|_\infty |\overline{u}^{k}|_{D^2} |\nabla \overline{u}^{k+1}|_2,\\
J_2=& \int_{\mathbb{R}^2} -\partial_{x}^\zeta(\overline{u}^{k} \cdot \nabla u^{k-1})\cdot \partial_{x}^\zeta\overline{u}^{k+1} \text{d}x\\
\leq &C|\nabla\overline{u}^{k}|_2 |\nabla\overline{u}^{k+1}|_2 |\nabla u^{k-1}|_\infty +C|\overline{u}^{k}|_6 |\nabla \overline{u}^{k+1}|_2 |\nabla^2 u^{k-1}|_3,\\
J_3=&\int_{\mathbb{R}^2} -\theta \partial_{x}^\zeta(\nabla ((\phi^{k+1})^2-(\phi^k)^2) )\cdot \partial_{x}^\zeta\overline{u}^{k+1} \text{d}x \\
 \leq &C |(\nabla\phi^{k+1}+\nabla\phi^{k})|_\infty|\nabla^2\overline{u}^{k+1}|_2 |\overline{\phi}^{k+1}|_2+C |(\phi^{k+1}+\phi^{k})|_\infty|\nabla^2\overline{u}^{k+1}|_2 |\nabla\overline{\phi}^{k+1}|_2,\\
J_4=&\int_{\mathbb{R}^2} \partial_{x}^\zeta\big(\psi^{k+1}\cdot Q(\overline{u}^k)\big)\cdot \partial_{x}^\zeta\overline{u}^{k+1}  \text{d}x\\
\leq& C|\psi^{k+1}|_\infty |\partial_{x}^2\overline{u}^{k}|_{2} |\nabla \overline{u}^{k+1}|_2+C|\nabla\psi^{k+1}|_3 |\nabla \overline{u}^{k}|_{6} |\nabla \overline{u}^{k+1}|_2,\\
J_5=&\int_{\mathbb{R}^2} \partial_{x}^\zeta\big(\overline{\psi}^{k+1}\cdot Q(u^{k-1})\big)\cdot \partial_{x}^\zeta\overline{u}^{k+1} \text{d}x
 \leq C|\overline{\psi}^{k+1}|_2 |\nabla u^{k-1}|_\infty |\nabla \partial_{x}^\zeta\overline{u}^{k+1}|_2.
\end{split}
\end{equation}
Using Young's inequality and (\ref{litong3}), we have
\begin{equation}\label{gogo12}\begin{split}
&\frac{d}{dt}| \nabla\overline{u}^{k+1}|^2_2+\alpha |\partial_{x}^2\overline{u}^{k+1} |^2_{2}\\
\leq& F^k_\eta(t)|\nabla\overline{u}^{k+1}|^2_2+F^k_2(t)\|\overline{\phi}^{k+1}\|^2_{1}
+F^k_3(t)|\overline{\psi}^{k+1}|^2_{2}+\eta\|\overline{u}^{k}\|^2_2,
\end{split}
\end{equation}
where
\begin{equation*}
\begin{split}
F^k_\eta(t)=&C\left(1+\frac{1}{\eta}\|u^{k}\|^2_{3}+\frac{1}{\eta}\|\nabla u^{k-1}\|^2_{2}+\frac{1}{\eta}|\psi^{k+1}|^2_{L^6\cap D^1\cap D^2}\right), \\ F^k_2(t)=&C(\|\phi^{k+1}\|_3+\|\phi^{k}\|_3)^2,\quad
F^k_3(t)=C\|\nabla u^{k-1}\|^2_2 ,
\end{split}
\end{equation*}
and
$$\int_0^t \big(F^k_\eta(s)+F^k_2(s)+F^k_3(s)\big)\text{d}s\leq C+C_\eta t,$$ for $t\in (0,T_{**}].$

Combining (\ref{gogo1}) and (\ref{gogo12}), we easily get
\begin{equation}\label{gogo13}\begin{split}
&\frac{d}{dt}\|\overline{u}^{k+1}\|^2_1+\alpha \|\nabla \overline{u}^{k+1} \|^2_{1}\\
\leq& \Theta^k_\eta(t)\|\overline{u}^{k+1}\|^2_1+\Theta^k_2(t)\|\overline{\phi}^{k+1}\|^2_{1}
+\Theta^k_3(t)|\overline{\psi}^{k+1}|^2_{2}+\eta\|\overline{u}^{k}\|^2_2,
\end{split}
\end{equation}
and
$$\int_0^t \big(\Theta^k_\eta(s)+\Theta^k_2(s)+\Theta^k_3(s)\big)\text{d}s\leq C+C_\eta t,$$
for $t\in (0,T_{**}]$.

Finally, let
\begin{equation*}\begin{split}
\Gamma^{k+1}(t)=&\sup_{0\leq s \leq t}\|\overline{\phi}^{k+1}(s)\|^2_{1}+\sup_{0\leq s \leq t}|\overline{\psi}^{k+1}(s)|^2_{ 2}+\sup_{0\leq s \leq t}\|\overline{u}^{k+1}(s)\|^2_1.
\end{split}
\end{equation*}
According to (\ref{go64qqqqqq}), (\ref{go64aa2}), (\ref{gogo13}) and  Gronwall's inequality, we have

\begin{equation*}\begin{split}
&\Gamma^{k+1}(t)+\int_{0}^{t}\mu\|\nabla\overline{u}^{k+1}\|^2_1\text{d}s
\leq  \Big( C\eta\int_{0}^{t}  \|\nabla \overline{u}^k\|^2_1\text{d}s+C\eta t\sup_{0\leq s \leq t}|\overline{u}^{k}(s)|^2_2 \Big)\exp{(C+C_\eta t)}.
\end{split}
\end{equation*}

We choose $\eta>0$ and $T_* \in (0,\min(1,T_{**}))$ small enough such that
$$
C\eta\exp{C}\leq \min\left(\frac{1}{4}, \frac{\mu}{4}\right), \quad \text{and}\quad \text{exp}(C_\eta T_*) \leq 2.
$$
Then we easily have
\begin{equation*}\begin{split}
\sum_{k=1}^{\infty}\Big(  \Gamma^{k+1}(T_*)+\int_{0}^{T_*} \mu\|\nabla\overline{u}^{k+1}\|^2_1\text{d}s\Big)\leq C<+\infty.
\end{split}
\end{equation*}
Thanks to
\begin{equation*}\begin{split}
\lim_{k\mapsto +\infty} |\overline{\psi}^{k+1}|_{6}\leq C\lim_{k\mapsto +\infty}\big( |\overline{\psi}^{k+1}|^{\frac{2}{3}}_{\infty} |\overline{\psi}^{k+1}|^{\frac{1}{3}}_{2}\big)\leq C\lim_{k\mapsto +\infty} |\overline{\psi}^{k+1}|^{\frac{1}{3}}_{2}=0,
\end{split}
\end{equation*}
we easily  know that the whole sequence $(\phi^k,\psi^k,u^k)$ converges to a limit $(\phi,\psi, u)$ in the following strong sense:
\begin{equation}\label{str}
\begin{split}
&\phi^k\rightarrow \phi\ \text{in}\ L^\infty([0,T_*];H^1(\mathbb{R}^2)),\\
&\psi^k\rightarrow\psi \ \text{in}\ L^\infty([0,T_*];L^6(\mathbb{R}^2)),\\
&u^k\rightarrow u\ \text{in}\ L^\infty ([0,T_*];H^1(\mathbb{R}^2)) \cap L^2([0,T_*];D^2(\mathbb{R}^2)).
\end{split}
\end{equation}

It is clear that $(\phi,\psi, u)$ satisfies the estimates in (\ref{jkk}).
Thanks to (\ref{str}), $(\phi,u)$ is a weak solution of problem \eqref{eq:cccq}-\eqref{qwe} with the following regularities:
\begin{equation}\label{rjkqq}\begin{split}
& \phi \in L^\infty([0,T_*];H^3),\quad  \phi_t \in L^\infty([0,T_*];H^2),\ \psi \in L^\infty([0,T_*] ; L^6\cap D^1\cap D^2),\\
&\partial_i \psi^{(j)}=\partial_j \psi^{(i)} \ (i,j=1,2), \ \psi_t \in L^\infty([0,T_*]; H^1),\ \psi_{tt} \in L^2([0,T_*]; L^2),\\
& u\in L^\infty([0,T_*]; H^3)\cap L^2([0,T_*] ; H^4), \ u_t \in L^\infty([0,T_*]; H^1)\cap L^2([0,T_*] ; D^2),\\
& u_{tt}\in L^2([0,T_*];L^2),\quad  t^{\frac{1}{2}}u\in L^\infty([0,T_*];D^4),\\
&t^{\frac{1}{2}}u_t\in L^\infty([0,T_*];D^2)\cap L^2([0,T_*];D^3) ,\ t^{\frac{1}{2}}u_{tt}\in L^\infty([0,T_*];L^2)\cap L^2([0,T_*];D^1).
\end{split}
\end{equation}
The time-continuity of the above solution can be obtained by standard procedure(see \cite{CK}). Therefore, it is a regular solution.

Now we prove the uniqueness of regular solutions. Let $(\phi_1,u_1)$ and $(\phi_2,u_2)$ be two regular solutions to  Cauchy problem (\ref{eq:cccq})-(\ref{qwe})  satisfying the uniform estimates in (\ref{jkk}). We denote that
$$
\overline{\phi}=\phi_1-\phi_2,\quad \overline{u}=u_1-u_2,
$$
and
$$\quad \overline{\psi}=\psi_1-\psi_2=\frac{2}{\gamma-1}\left(\nabla \phi_{1}/\phi_{1}-\nabla \phi_{2}/\phi_{2}\right).$$
Then $(\overline{\phi},\overline{\psi},\overline{u})$ satisfies the system
 \begin{equation}
\label{zhuzhu}
\begin{cases}
\ \ \displaystyle \overline{\phi}_t+u_1\cdot \nabla\overline{\phi} +\overline{u}\cdot\nabla\phi_{2}+\frac{\gamma-1}{2}(\overline{\phi} \text{div}u_2 +\phi_{1}\text{div}\overline{u})=0,\\[8pt]
\ \ \displaystyle \overline{\psi}_t+\sum_{l=1}^2 A_l(u_1) \partial_l\overline{\psi}+B(u_{1})\overline{\psi}+\nabla \text{div}\overline{u}=\overline{\Upsilon}_1+\overline{\Upsilon}_2,\\[8pt]
\ \ \overline{u}_t+ u_1\cdot\nabla \overline{u}+ \overline{u}\cdot \nabla u_{2}+\theta\nabla ((\phi_1)^2-(\phi_2)^2) \\[8pt]
\qquad=-L\overline{u}+\psi_1\cdot Q(\overline{u})+\overline{\psi}\cdot Q(u_2),
\end{cases}
\end{equation}
with $\overline{\Upsilon}_1$  and $\overline{\Upsilon}_2$ defined by
\begin{equation*}
\overline{\Upsilon}_1=-\sum_{l=1}^2(A_l(u_{1}) \partial_l\psi_{2}-A_l(u_{2}) \partial_l\psi_{2}),\quad \overline{\Upsilon}_2=-(B(u_{1}) \psi_{2}-B(u_{2}) \psi_{2}).
\end{equation*}
Let
$$
\Phi(t)=\|\overline{\phi}(t)\|^2_{1}+|\overline{\psi}(t)|^2_{ 2}+\|\overline{u}(t)\|^2_1.
$$
Similarly to the derivation of (\ref{go64})-(\ref{gogo1}), we can show that
\begin{equation}\label{gonm}\begin{split}
\frac{d}{dt}\Phi(t)+C\|\nabla \overline{u}(t)\|^2_1\leq G(t)\Phi (t),
\end{split}
\end{equation}
where $ \int_{0}^{t}G(s)ds\leq C$, for $0\leq t\leq T_*$. From the Gronwall's inequality, we conclude that
$$\overline{\phi}=\overline{\psi}=\overline{u}=0,$$
then the uniqueness is obtained.

\subsection{Proof of Theorem \ref{th2} and Corollary \ref{co2}}\label{3.6}\ \\

Based on Theorem \ref{th1}, we are now ready to prove the local existence of regular solution to the original Cauchy problem  (\ref{eq:1.1}). Moreover, we will show that the regular solutions that we obtained satisfy system (\ref{eq:1.1}) classically.\\

\textbf{Proof of Theorem \ref{th2}.}

\begin{proof}
For initial data (\ref{th78}),  we know from Theorem \ref{th1} that there exists  a time $T_{*}> 0$ such that the problem (\ref{eq:cccq})-(\ref{qwe}) has a unique regular solution $(\phi,u)$ satisfying the regularity (\ref{reg11qq}), which means that
\begin{equation}\label{reg2}
\begin{split}
(\rho^{\frac{\gamma-1}{2}},u )=(\phi,u)\in C^1([0,T_{*}]\times \mathbb{R}^2), \quad (\nabla \rho/\rho,\partial_{x}^\xi u) \in   C([0,T_{*}]\times \mathbb{R}^2),
\end{split}
\end{equation}
where $\xi\in \mathbb{R}^2$ with $|\xi|=2$. Since
$$
\rho(t,x)=\phi^{\frac{2}{\gamma-1}}(t,x),
$$
and $ \frac{2}{\gamma-1}\geq 1$ for $ 1< \gamma \leq3 $, it is easy to show that
$$\rho(t,x)\in C^1([0,T_{*}]\times\mathbb{R}^2).$$

Multiplying $(\ref{eq:cccq})_1$ by
$
\frac{\partial \rho}{\partial \phi}(t,x)=\frac{2}{\gamma-1}\phi^{\frac{3-\gamma}{\gamma-1}}(t,x)\in C([0,T_{*}]\times \mathbb{R}^2)
$,
we get the continuity equation in (\ref{eq:1.1}):
\begin{equation} \label{eq:2.58}
\rho_t+u \cdot\nabla \rho+\rho\text{div} u=0.
\end{equation}

Multiplying $(\ref{eq:cccq})_2$ by
$
\phi^{\frac{2}{\gamma-1}}=\rho(t,x)\in C^1([0,T_{*}]\times \mathbb{R}^2)
$,
we get the momentum equations in (\ref{eq:1.1}):
\begin{equation} \label{eq:2.60}
\begin{split}
&\rho u_t+\rho u\cdot \nabla u+\nabla P=\text{div}\Big(\mu(\rho)(\nabla u+ (\nabla u)^\top)+\lambda(\rho)\text{div}u I_2\Big).
\end{split}
\end{equation}
That is to say, $(\rho,u)$ satisfies problem  (\ref{eq:1.1})  in classical sense with regularity (\ref{reg11}).

Recalling that $\rho$ can be represented by the formula
$$
\rho(t,x)=\rho_0(W(0,t,x))\exp\big(\int_{0}^{t}\textrm{div} u(s,W(s,t,x))\text{d}s\big),
$$
where  $W\in C^1([0,T_{*}]\times[0,T_{*}]\times \mathbb{R}^2)$ is the solution to the initial value problem
\begin{equation}
\label{eq:bb1}
\begin{cases}
\frac{d}{dt}W(t,s,x)=u(t,W(t,s,x)),\quad 0\leq t\leq T_{*},\\[4pt]
W(s,s,x)=x, \quad\quad\quad 0\leq s\leq T_{*},\quad  x\in \mathbb{R}^2,
\end{cases}
\end{equation}
it is obvious that
$$
\rho(t,x)\geq 0, \ \forall (t,x)\in [0,T_{*}]\times \mathbb{R}^2.
$$
In summary, the Cauchy problem  (\ref{eq:1.1}) has a unique regular solution $(\rho,u)$.
\end{proof}

\textbf{Proof of Corollary \ref{co2}}.
\begin{proof} When $1< \gamma \leq \frac{5}{3}$, $\frac{2}{\gamma-1}\geq 3$. Since
$\phi \in C([0,T_*],H^{3})\bigcap C^1([0,T_*],H^{2})$, we read
$$
\rho(t,x)=\phi^{\frac{2}{\gamma-1}}(t,x),
$$
that
$$
\rho \in C([0,T_*],H^{3}).
$$
With the help of the continuity equation
\begin{equation*}
\rho_t+u \cdot\nabla \rho+\rho\text{div} u=0,
\end{equation*}
and the fact that
$
u(t,x)\in C([0,T_*],H^{3})\bigcap C^1([0,T_*],H^{2})$,  it is clear that $$
\rho \in C([0,T_*],H^{3})\cap C^1([0,T_*],H^{2}),
$$
and the regularity on $\rho_{t}$ in Corollary \ref{co2} follows.

Furthermore, when  $\gamma=2$, or $3$, $\rho(t,x)=\phi^2(t,x)$ and $\rho(t,x)=\phi(t,x)$,  respectively. By the same token, the regularity of $\rho$ in these cases can be achieved. \end{proof}

\section{Stability in $H^2$ sense}
Now we prove the stability  in $H^2$, i.e., Theorem \ref{co02}. For $i=1,2$, let $(\phi_i,u_i)$  be the regular solution to  Cauchy problem (\ref{eq:cccq}) with initial data $(\phi_{0i},u_{0i})$ satisfying \eqref{th78}.  Let $K>0$ be a constant such that
$$
\|\phi_{0i}\|_2+\left|\frac{\nabla \phi_{0i}}{\phi_{0i}}\right|_{L^6\cap D^1}+\|u_{0i}\|_2\leq K.
$$
Denote
$$
\tilde{\phi}=\phi_1-\phi_2,\quad \tilde{u}=u_1-u_2,
$$
and
$$\quad \tilde{\psi}=\psi_1-\psi_2=\frac{2}{\gamma-1}\left(\frac{\nabla \phi_{1}}{\phi_{1}}-\frac{\nabla \phi_{2}}{\phi_{2}}\right).$$
Then $(\tilde{\phi},\tilde{\psi},\tilde{u})$ satisfies the system
 \begin{equation}
\label{zhuzhu}
\begin{cases}
\ \ \displaystyle \tilde{\phi}_t+u_1\cdot \nabla\tilde{\phi} +\tilde{u}\cdot\nabla\phi_{2}+\frac{\gamma-1}{2}(\tilde{\phi} \text{div}u_2 +\phi_{1}\text{div}\tilde{u})=0,\\[8pt]
\ \ \displaystyle \tilde{\psi}_t+\sum_{l=1}^2 A_l(u_1) \partial_l\tilde{\psi}+B(u_{1})\tilde{\psi}+\nabla \text{div}\tilde{u}=\tilde{\Upsilon}_1+\tilde{\Upsilon}_2,\\[8pt]
\ \ \tilde{u}_t+ u_1\cdot\nabla \tilde{u}+ \tilde{u}\cdot \nabla u_{2}+\theta\nabla ((\phi_1)^2-(\phi_2)^2) \\[8pt]
\qquad=-L\tilde{u}+\psi_1\cdot Q(\tilde{u})+\tilde{\psi}\cdot Q(u_2),
\end{cases}
\end{equation}
where $\tilde{\Upsilon}_1$  and $\tilde{\Upsilon}_2$ are defined by
\begin{equation*}
\tilde{\Upsilon}_1=-\sum_{l=1}^2(A_l(u_{1}) \partial_l\psi_{2}-A_l(u_{2}) \partial_l\psi_{2}),\quad \tilde{\Upsilon}_2=-(B(u_{1}) \psi_{2}-B(u_{2}) \psi_{2}).
\end{equation*}

Similarly to the derivation of (\ref{go64qqqqqq}), we have
\begin{equation}\label{go99gg}
\begin{cases}
\displaystyle
\frac{d}{dt}\| \tilde{\phi}(t)\|^2_1\leq A_\eta(t)\| \tilde{\phi}(t)\|^2_1+C \|\tilde{u}(t)\|^2_2+\eta \|\nabla\tilde{u}(t)\|^2_1,\\[8pt]
\displaystyle
 \int_0^t A_\eta(s)\text{d}s\leq C+C_{\eta}t\quad \text{for}\quad  t\in[0,T_*].
\end{cases}
\end{equation}
Then taking $\partial_{x}^\xi$ to $(\ref{zhuzhu})_1$ ($|\xi|=2$),
multiplying by $2 \partial_{x}^\xi\tilde{\phi}$ and integrating over $\mathbb{R}^2$,  we have
\begin{equation}\label{ghu}
\begin{split}
\frac{d}{dt}|\partial_{x}^\xi\tilde{\phi}|^2_{2}=&-2 \int_{\mathbb{R}^2}\partial_{x}^\xi\Big(u_1\cdot \nabla\tilde{\phi} +\tilde{u}\cdot\nabla\phi_2 +\frac{\gamma-1}{2}(\tilde{\phi} \text{div}u_2+\phi_1 \text{div}\tilde{u})\Big)\partial_{x}^\xi\tilde{\phi}\text{d}x\\
\leq & C(\|\nabla u_1\|_2+\|\nabla u_2\|_2+\| \phi_1\|^2_3+\|\nabla \phi_2\|^2_2)\|\tilde{\phi}\|^2_{2}+C\|\tilde{u}\|^2_{2}+\eta \|\nabla \tilde{u}\|^2_2.
\end{split}
\end{equation}

(\ref{go99gg})-(\ref{ghu}) yield
\begin{equation}\label{go64ppppp}
\begin{cases}
\displaystyle
\frac{d}{dt}\| \tilde{\phi}(t)\|^2_2\leq \Phi_\eta(t)\| \tilde{\phi}(t)\|^2_2+C\|\tilde{u}(t)\|^2_{2}+\eta \|\nabla\tilde{u}(t)\|^2_2,\\[8pt]
\displaystyle
 \int_0^t \Phi_\eta(s)\text{d}s\leq C+C_{\eta}t\quad \text{for}\quad  t\in[0,T_*].
\end{cases}
\end{equation}

Multiplying $(\ref{zhuzhu})_1$ by $6|\tilde{\phi}|^4\tilde{\phi}$ and integrating over $\mathbb{R}^2$,
\begin{equation}\label{goklpo}
\begin{split}
\frac{d}{dt}|\tilde{\psi}|^6_6\leq& C\Big(\sum_{l=1}^{2}|\partial_{l}A_l(u_{1})|_\infty+|B(u_{1})|_\infty\Big)|\tilde{\psi}|^6_6+|\nabla^2 \tilde{u}|_6|\tilde{\psi}|^5_6+C(|\tilde{\Upsilon}_1 |_6+|\tilde{\Upsilon}_{2}|_6)|\tilde{\psi}|^5_6.
\end{split}
\end{equation}
From H\"older's inequality, we easily have
\begin{equation}\label{zhuuu}
|\tilde{\Upsilon}_1 |_6+|\tilde{\Upsilon}_{2}|_6\leq C(|\tilde{u}|_\infty |\nabla \psi_2|_6+|\psi_2|_\infty |\nabla \tilde{u}|_6)\leq C|\psi_2|_{L^6\cap D^1\cap D^2}\|\tilde{u}\|_2,
\end{equation}
and
\begin{equation}\label{gokluu}
\begin{cases}
\displaystyle
\frac{d}{dt}|\tilde{\psi}|^2_6\leq B_\eta(t) |\tilde{\psi}(t)|^2_6+C|\tilde{u}(t)|^2_2+\eta\|\nabla\tilde{u}(t) \|^2_2,\\[8pt]
\displaystyle
\int_0^t B_\eta(s)\text{d}s\leq C+C_{\eta}t \quad \text{for} \quad t\in[0,T_*].
\end{cases}
\end{equation}

Taking $\partial_{x}^\zeta$ to $(\ref{zhuzhu})_2$ ($|\zeta|=1$), we have
\begin{equation}\label{str11}
\begin{split}
&\partial_{x}^\zeta\tilde{\psi}_t+\sum_{l=1}^2 A_l(u_{1}) \partial_l \partial_{x}^\zeta\tilde{\psi}+B(u_{1})\partial_{x}^\zeta\tilde{\psi}+\partial_{x}^\zeta\nabla \text{div}\tilde{u}^k\\
=&\partial_{x}^\zeta\tilde{\Upsilon}_1+\partial_{x}^\zeta\tilde{\Upsilon}_2+\Upsilon_{11}+\Upsilon_{22},
\end{split}
\end{equation}
where
\begin{equation*}
\Upsilon_{11}=\sum_{l=1}^2 \big(A_l(u_{1}) \partial_l \partial_{x}^\zeta\tilde{\psi}-\partial_{x}^\zeta(A_l(u_{1}) \partial_l \tilde{\psi})\big),\
\Upsilon_{22}=B(u_{1})\partial_{x}^\zeta\tilde{\psi}-\partial_{x}^\zeta(B(u_{1})\tilde{\psi}).
\end{equation*}
 Multiplying  (\ref{str11}) by $2\partial_{x}^\zeta\tilde{\psi}$ and integrating over $\mathbb{R}^2$, we have
\begin{equation}\label{go64aa66}
\begin{split}
\frac{d}{dt}|\partial_{x}^\zeta\tilde{\psi}|^2_2\leq& C\Big(\sum_{l=1}^{2}|\partial_{l}A_l(u_{1})|_\infty+|B(u_{1})|_\infty\Big)|\partial_{x}^\zeta\tilde{\psi}|^2_2+|\nabla^3 \tilde{u}|_2|\partial_{x}^\zeta\tilde{\psi}|_2\\
&+C(|\partial_{x}^\zeta\tilde{\Upsilon}_1 |_2+|\partial_{x}^\zeta\tilde{\Upsilon}_2|_2+|\Upsilon_{11} |_2+|\Upsilon_{22}|_2)|\partial_{x}^\zeta\tilde{\psi}|_2.
\end{split}
\end{equation}
From H\"older's inequality and Lemma \ref{zhen1}, it is easy to deduce that
\begin{equation}\label{go64aa166}
|\partial_{x}^\zeta\tilde{\Upsilon}_1 |_2+| \partial_{x}^\zeta\tilde{\Upsilon}_2 |_2\leq C|\psi_{2}|_{L^6\cap D^1 \cap D^2} \|\nabla \tilde{u}\|_1,\
| \Upsilon_{11} |_2+| \Upsilon_{22} |_2\leq C|\nabla \tilde{\psi}|_2 \|\nabla u_{1}\|_3.
\end{equation}
From (\ref{gokluu})-(\ref{go64aa166}), we have
\begin{equation}\label{go64ll2}\begin{cases}
\displaystyle
\frac{d}{dt}|\tilde{\psi}(t)|^2_{L^6 \cap D^1}\leq \Psi_\eta(t)|\tilde{\psi}(t)|^2_{L^6\cap D^1}+C\|\tilde{u}(t)\|^2_2+\eta \|\nabla\tilde{u}(t)\|^2_2,\\[8pt]
\displaystyle
 \int_0^t \Psi_\eta(s)\text{d}s\leq C+C_{\eta}t \quad \text{for} \quad t\in[0,T_*].
\end{cases}
\end{equation}

Similarly to the derivation of (\ref{gogo13}),  we easily have
\begin{equation}\label{gogo16y}\begin{split}
\displaystyle
&\frac{d}{dt}\|\tilde{u}(t)\|^2_1+ \|\nabla \tilde{u}(t) \|^2_{1}\\
\leq&E_\eta(t)\|\tilde{u}(t)\|^2_1+E_2(t)\|\tilde{\phi}(t)\|^2_{1}
+E_3(t)|\tilde{\psi}(t)|^2_{L^6\cap D^1}+\eta \|\nabla\tilde{u}(t)\|^2_2,
\end{split}
\end{equation}
where we have $
\int_0^t \big(E_\eta(s)+E_2(s)+E_3(s)\big)\text{d}s\leq C+C_{\eta}t \ \text{for}\  t\in (0,T_*]
$.

Taking $\partial_{x}^{\xi}$ to $(\ref{zhuzhu})_3$ ($|\xi|=2$), multiplying by $\partial_{x}^\xi\tilde{u}$ and  integrating over $\mathbb{R}^2$,  we have
\begin{equation*}
\begin{split}
&\frac{1}{2}\frac{d}{dt}|\partial_{x}^\xi\tilde{u}|^2_{2}+\alpha|\nabla \partial_{x}^\xi\tilde{u} |^2_2+(\alpha+\beta)| \partial_{x}^\xi \text{div} \tilde{u} |^2_2\\
=& \int_{\mathbb{R}^2} \Big(-\partial_{x}^\xi( u_1\cdot\nabla \tilde{u})-\partial_{x}^\xi(\tilde{u} \cdot \nabla u_2)-\theta \partial_{x}^\xi(\nabla (\phi^2_1-\phi^2_2) )\Big)\cdot \partial_{x}^\xi\tilde{u} \text{d}x \\
&+\int_{\mathbb{R}^2} \Big(\partial_{x}^\xi\big(\psi_1\cdot Q(\tilde{u})\big)+\partial_{x}^\xi\big(\tilde{\psi}\cdot Q(u_2)\big)\Big)\cdot \partial_{x}^\xi\tilde{u} \text{d}x=\sum_{i=6}^{10} J_i.
\end{split}
\end{equation*}
The right-hand side can be estimated term by term as follows.
\begin{equation}
\begin{split}\label{litongty}
J_6=&\int_{\mathbb{R}^2} -\partial_{x}^\xi( u_1\cdot\nabla \tilde{u})\cdot \partial_{x}^\xi\tilde{u} \text{d}x\\
\leq &C|\nabla^2 u_1|_3|\nabla \tilde{u}|_6 |\nabla^2\tilde{u}|_2+C| u_1|_\infty |\nabla^3\tilde{u}|_{2} |\nabla^2 \tilde{u}|_2+C| \nabla u_1|_\infty |\nabla^2\tilde{u}|^2_{2},\\
J_7=& \int_{\mathbb{R}^2} -\partial_{x}^\xi(\tilde{u} \cdot \nabla u_2)\cdot \partial_{x}^\xi\tilde{u} \text{d}x\\
\leq &C|\nabla u_2|_\infty |\nabla^2\tilde{u}|^2_2 +C |\nabla^2 u_2|_3|\nabla \tilde{u}|_6|\nabla^2 \tilde{u}|_2+C |\nabla^3 u_2|_3| \tilde{u}|_6|\nabla^2 \tilde{u}|_2,\\
J_8=&\int_{\mathbb{R}^2} -\theta \partial_{x}^\xi(\nabla (\phi^2_1-\phi^2_2) )\cdot \partial_{x}^\xi\tilde{u} \text{d}x \\
 \leq &C\big( |\nabla(\phi_1+\phi_2)|_\infty|\nabla\tilde{\phi}|_2+ |(\phi_1+\phi_2)|_\infty |\nabla^2\tilde{\phi}|_2+ |\nabla^2(\phi_1+\phi_2)|_3 |\tilde{\phi}|_6\big)|\nabla^3\tilde{u}|_2,\\
J_9=&\int_{\mathbb{R}^2} \partial_{x}^\xi\big(\psi_1\cdot Q(\tilde{u})\big)\cdot \partial_{x}^\xi\tilde{u}  \text{d}x\\
\leq& C|\psi_1|_\infty |\nabla^3\tilde{u}|_{2} |\nabla^2 \tilde{u}|_2+C|\nabla\psi_1|_3 |\nabla^2\tilde{u}|_{6} |\nabla ^2 \tilde{u}|_2+C|\nabla^2\psi_1|_2 |\nabla\tilde{u}|_{3} |\nabla ^2 \tilde{u}|_6,\\
J_{10}=&\int_{\mathbb{R}^2} \partial_{x}^\xi\big(\tilde{\psi}\cdot Q(u_2)\big)\cdot \partial_{x}^\xi\tilde{u} \text{d}x
 \leq C|\nabla\tilde{\psi}|_2 |\nabla u_2|_\infty |\nabla^3\tilde{u}|_2+C|\tilde{\psi}|_6 |\nabla^2 u_2|_3 |\nabla^3\tilde{u}|_2.
\end{split}
\end{equation}
According to Young's inequality and (\ref{litongty}), we have
\begin{equation}\label{gogo12y}\begin{cases}
\displaystyle
\frac{d}{dt}| \tilde{u}(t)|^2_{D^2}+\alpha |\tilde{u} (t)|^2_{D^3}
\leq F_\eta(t)\|\tilde{u}(t)\|^2_{D^2}+F_2(t)\|\tilde{\phi}(t)\|^2_{2}
+F_3(t)|\tilde{\psi}(t)|^2_{L^6 \cap D^1},\\[10pt]
\displaystyle
\int_0^t \big(F_\eta(s)+F_2(s)+F_3(s)\big)\text{d}s\leq C+C_{\eta}t\quad \text{for}\quad t\in (0,T_*].
\end{cases}
\end{equation}
Then combining (\ref{gogo16y}) and (\ref{gogo12y}), we easily have
\begin{equation}\label{gogo19}\begin{cases}
\displaystyle
\frac{d}{dt}\|\tilde{u}(t)\|^2_2+\alpha \|\nabla \tilde{u}(t) \|^2_{2}
\leq \Theta_\eta(t)\|\tilde{u}(t)\|^2_2+\Theta_2(t)\|\tilde{\phi}(t)\|^2_{2}
+\Theta_3(t)|\tilde{\psi}(t)|^2_{L^6\cap D^1},\\[10pt]
\displaystyle
\int_0^t \big(\Theta_\eta(s)+\Theta_2(s)+\Theta_3(s)\big)\text{d}s\leq C+C_{\eta}t\quad \text{for}\quad t\in (0,T_*].
\end{cases}
\end{equation}

Finally, let
\begin{equation*}\begin{split}
\Gamma(t)=&\|\tilde{\phi}(t)\|^2_{2}+|\tilde{\psi}(t)|^2_{ L^6\cap D^1}+\|\tilde{u}(t)\|^2_2.
\end{split}
\end{equation*}
Then we have
\begin{equation*}\begin{split}
&\frac{d}{dt}\Gamma(t)+\mu\|\nabla\tilde{u}(t) \|^2_2
\leq \Pi'_\eta(t) \Gamma(t)
\end{split}
\end{equation*}
for some $\Pi'_\eta$ such that  $\int_{0}^{t}\Pi'_\eta(s)\text{d}s\leq C+C_{\eta}t$ for $ t\in (0,T_*]$.
Then our stability result follows from Gronwall's inequality.

\section{Blow-up criterion}

In this section, we give the proof to Theorem \ref{th3s}. In order to prove (\ref{cri}), we use a contradiction argument. Let $( \rho, u)$ be the unique regular solution to  the Cauchy problem  (\ref{eq:1.1}) with the maximal existence time $\overline{T}$. We assume that
$\overline{T}<+\infty$ and
\begin{equation}\label{we11}
\begin{split}
\lim_{T\mapsto \overline{T}} \left(\sup_{0\leq t \leq T}\left|\frac{\nabla \rho} { \rho} (t)\right|_{6}+\int_{0}^{T}|D(u(t))|_{\infty}\text{d}t\right)=C_0<+\infty.
\end{split}
\end{equation}
We will show that under assumption \eqref{we11}, $\overline{T}$ is actually not the maximal existence time for the regular solution.

From the definition of regular solutions, we know that, for $\phi=\rho^{\frac{\gamma-1}{2}}$,
$(\phi, u)$ satisfies
\begin{equation}
\begin{cases}
\label{we22}
\phi_t+u\cdot \nabla \phi+\frac{\gamma-1}{2}\phi\text{div} u=0,\\[12pt]
\psi_t+\nabla (u\cdot \psi)+\nabla \text{div} u=0,\\[12pt]
u_t+u\cdot\nabla u +2\theta\phi\nabla \phi+Lu=\psi\cdot Q(u).
 \end{cases}
\end{equation}
For $(\ref{we22})_2$, we have the equivalent form
\begin{equation}\label{kucccc}
\psi_t+\sum_{l=1}^2 A_l \partial_l\psi+B\psi+\nabla \text{div} u=0.
\end{equation}
Here $A_l=(a^{(l)}_{ij})_{2\times 2}$ ($i,j,l=1,2$) are symmetric  with $a^{(l)}_{ij}=u^{(l)}$ when $i=j$; and $a^{(l)}_{ij}=0$, otherwise. $B=(\nabla u)^\top$, so (\ref{kucccc}) is a positive symmetric hyperbolic system.

By assumptions (\ref{we11}) and (\ref{we22}), we first show that  density $\rho$ is uniformly bounded.
 \begin{lemma}\label{s1} Let $( \rho, u)$ be the unique regular solution to  the Cauchy problem  (\ref{eq:1.1}) on $[0,\overline{T})$ satisfying \eqref{we11}. Then
\begin{equation*}
\begin{split}
\|\rho\|_{L^\infty([0,T]\times \mathbb{R}^2)}+\|\phi\|_{L^\infty([0,T]; L^q( \mathbb{R}^2))}\leq C, \quad 0\leq T< \overline{T},
\end{split}
\end{equation*}
where $C>0$ depends on $C_0$,  constant $q\in [2,+\infty]$ and $T$.
 \end{lemma}
 \begin{proof}
First, it is obvious that  $\phi$ can be represented by
\begin{equation}
\label{eq:bb1a}
\phi(t,x)=\phi_0(W(0,t,x))\exp\Big(-\frac{\gamma-1}{2}\int_{0}^{t}\textrm{div}u(s,W(s,t,x))\text{d}s\Big),
\end{equation}
where  $W\in C^1([0,T]\times[0,T]\times \mathbb{R}^2)$ is the solution to the initial value problem
\begin{equation}
\label{eq:bb1z}
\begin{cases}
\frac{d}{dt}W(t,s,x)=u(t,W(t,s,x)),\quad 0\leq t\leq T,\\[6pt]
W(s,s,x)=x, \quad \ \quad \quad 0\leq s\leq T,\ x\in \mathbb{R}^2.
\end{cases}
\end{equation}
Then it is clear that $\|\phi\|_{L^\infty([0,T]\times \mathbb{R}^2)}\leq |\phi_0|_\infty \exp({CC_0})$.

Next,  multiplying $(\ref{we22})_1$ by $2\phi$ and integrating over $\mathbb{R}^2$, we get
\begin{equation}\label{zhumw}
\begin{split}
& \frac{d}{dt}|\phi|^2_2\leq C |\text{div} u|_\infty|\phi|^2_2,
\end{split}
\end{equation}
 from (\ref{we11}) and  Gronwall's inequality, we immediately obtain the desired conclusions.
 \end{proof}
 Now we give the basic energy estimates.
  \begin{lemma}\label{s2} Let $( \rho, u)$ be the unique regular solution to  the Cauchy problem  (\ref{eq:1.1}) on $[0,\overline{T})$ satisfying \eqref{we11}. Then
\begin{equation*}
\begin{split}
\sup_{0\leq t\leq T}|u(t)|^2_{ 2}+\int_{0}^{T}|\nabla u(t)|^2_{2}\text{d}t\leq C,\quad 0\leq T<  \overline{T},
\end{split}
\end{equation*}
where $C$ only depends on $C_0$ and $T$.
 \end{lemma}
\begin{proof}  Multiplying $(\ref{we22})_3$ by $2u$ and integrating over $\mathbb{R}^2$, we have
\begin{equation}\label{zhu1}
\begin{split}
& \frac{d}{dt}|u|^2_2+\alpha|\nabla u|^2_2+(\alpha+\beta)|\text{div} u|^2_2\\
=&\int_{\mathbb{R}^2} 2\Big(-u\cdot \nabla u \cdot u-\theta \nabla \phi^2\cdot u+\psi \cdot Q(u)\cdot u \Big) \text{d}x\equiv: L_1+L_2+L_3.
\end{split}
\end{equation}
The right-hand side terms can be estimated as follows.
 \begin{equation}\label{zhu2s}
\begin{split}
L_1=&-\int_{\mathbb{R}^2} 2u\cdot \nabla u \cdot u \text{d}x\leq C|\text{div} u|_\infty|u|^2_2,\\
L_2=&\int_{\mathbb{R}^2} \theta \phi^2 \text{div} u \text{d}x\leq C|\phi|_2^{2}|\text{div} u|_\infty
\le  C|\text{div} u|_\infty,\\
L_3=&\int_{\mathbb{R}^2} 2\psi \cdot Q(u)\cdot u  \text{d}x\leq C|\psi|_6|\nabla u|_2|u|_3\\
\leq& \frac{\alpha}{4}|\nabla u|^2_2+C|\psi|^2_6|u|^{\frac{4}{3}}_2|\nabla u|^{\frac{2}{3}}_2
\leq \frac{\alpha}{2}|\nabla u|^2_2+C|\psi|^3_6|u|^2_2,
\end{split}
\end{equation}
where we have used the fact $|u|_3 \leq C|u|^{\frac{2}{3}}_2|\nabla u|^{\frac{1}{3}}_2$.  \eqref{zhu1} and  (\ref{zhu2s}) yield
\begin{equation}\label{zhu3s}
\begin{split}
& \frac{d}{dt}|u|^2_2+\frac{\alpha}{2}|\nabla u|^2_2\leq C(|\text{div} u|_\infty+1)|u|^2_2+C|\text{div} u|_\infty.
\end{split}
\end{equation}
According to Gronwall's inequality, we have
\begin{equation}\label{zhu5}
\begin{split}
|u(t)|^2_2+\int_0^t|\nabla u(s)|^2_2\text{d}s\leq C,\quad 0\leq t\leq T.
\end{split}
\end{equation}
\end{proof}

The next lemma is a key estimate on $\nabla \phi$ and $\nabla u$. We denote $\omega=\partial_{x_1} u^{(2)}-\partial_{x_2} u^{(1)}$.
  \begin{lemma}\label{s4} Let $( \rho, u)$ be the unique regular solution to  the Cauchy problem  (\ref{eq:1.1}) on $[0,\overline{T})$ satisfying \eqref{we11}. Then
\begin{equation*}
\begin{split}
\sup_{0\leq t\leq T}|\nabla u(t)|^2_{ 2}+\sup_{0\leq t\leq T}|\nabla \phi(t)|^2_{ 2}+\int_0^T (|\nabla^2 u|^2_2+| u_t|^2_2)\text{d}t\leq C,\quad 0\leq T<  \overline{T},
\end{split}
\end{equation*}
where $C$ only depends on $C_0$ and $T$.
 \end{lemma}
\begin{proof}
First, multiplying $(\ref{we22})_3$ by $-Lu-\theta \nabla \phi^2 $ and integrating over $\mathbb{R}^2$, we have
\begin{equation}\label{zhu6}
\begin{split}
&\frac{1}{2} \frac{d}{dt}\Big(\alpha|\nabla u|^2_2+(\alpha+\beta)|\text{div}u|^2_2\Big)+\int_{\mathbb{R}^2}(-Lu-\theta \nabla \phi^2)^2\text{d}x\\
=&-\alpha\int_{\mathbb{R}^2} (u\cdot \nabla u) \cdot \omega'\text{d}x+(2\alpha+\beta)\int_{\mathbb{R}^2} (u\cdot \nabla u) \cdot  \nabla \text{div}u\text{d}x\\
&-\theta\int_{\mathbb{R}^2} (u\cdot \nabla u) \cdot \nabla \phi^2 \text{d}x-\theta\int_{\mathbb{R}^2} u_t \cdot \nabla \phi^2 \text{d}x+\theta\int_{\mathbb{R}^2}(\psi\cdot Q(u)) \cdot  \nabla \phi^2\text{d}x\\
&-\int_{\mathbb{R}^2} (\psi\cdot Q(u)) \cdot (\alpha \triangle u+(\alpha+\beta)\nabla \text{div}u)\text{d}x\equiv: \sum_{i=4}^{9} L_i,
\end{split}
\end{equation}
where we have used the fact that $-\triangle u+\nabla\text{div}u=(\partial_{x_2} \omega, -\partial_{x_1} \omega)^\top=\omega'$.

From the standard elliptic estimate shown in Lemma \ref{zhenok}, we have
\begin{equation}\label{gaibian}
\begin{split}
&|\nabla^2 u|^2_2-C|\theta\nabla\phi^2|^2_2\\
\leq& C|\alpha \triangle u+(\alpha+\beta)\nabla \text{div}u|^2_2-C|\theta\nabla\phi^2|^2_2\\
\leq& C|\alpha \triangle u+(\alpha+\beta)\nabla \text{div}u-\theta\nabla\phi^2|^2_2.
\end{split}
\end{equation}

Now we estimate the right-hand side of (\ref{zhu6}) term by term.
According to
\begin{equation*}
\begin{split}
\frac{1}{2}\nabla (| u|^2)-u \cdot \nabla u=(u^{(2)}\omega,-u^{(1)}\omega)^\top=\omega'',
\end{split}
\end{equation*}
and H\"older's inequality, Gagliardo-Nirenberg inequality and Young's inequality,
 we  obtain
\begin{equation}\label{zhu10t}
\begin{split}
|L_4|=&\alpha\Big|\int_{\mathbb{R}^2} (u \cdot \nabla u) \cdot \omega'\text{d}x\Big|
=\alpha\Big| \int_{\mathbb{R}^2} \big(\frac{1}{2}\nabla (| u|^2)-\omega''\big)\cdot \omega'\text{d}x\Big|\\
=&\alpha\Big| \int_{\mathbb{R}^2} -\omega'' \cdot \omega'\text{d}x\Big|=\frac{\alpha}{2}\Big| \int_{\mathbb{R}^2} u^{(2)} \partial_{x_2}\omega^2+u^{(1)}\partial_{x_1}\omega^2\text{d}x\Big|\\
=&\frac{\alpha}{2}\Big| \int_{\mathbb{R}^2} \omega^2\text{div}u \text{d}x\Big|
\leq C|\text{div}u|_\infty|\nabla u|^2_2,\\
|L_6|=&\theta\Big|\int_{\mathbb{R}^2} (u\cdot \nabla u) \cdot \nabla \phi^2 \text{d}x\Big|\\
=&\theta\Big|-\int_{\mathbb{R}^2} \nabla u: (\nabla u)^\top \phi^2 \text{d}x-\int_{\mathbb{R}^2} \phi^2 u \cdot \nabla (\text{div}u)  \text{d}x\Big|\\
=&\theta\Big|-\int_{\mathbb{R}^2} \nabla u: (\nabla u)^\top\phi^2 \text{d}x+\int_{\mathbb{R}^2} (\text{div}u)^2 \phi^2 \text{d}x+\int_{\mathbb{R}^2}  u\cdot\nabla\phi^2\text{div}u \text{d}x\Big|\\
\leq &C|\nabla u|^2_2+C|\text{div}u|_{\infty}|u|_2|\nabla\phi^{2}|_2\\
\leq& C(|\nabla u|^2_2+|\text{div}u|_{\infty}+|\text{div}u|_{\infty}|\nabla\phi|_2^{2}),\\
\end{split}
\end{equation}
\begin{equation}\label{zhu10t33}
\begin{split}
L_{7}=&-\theta\int_{\mathbb{R}^2} u_t \cdot \nabla \phi^2 \text{d}x=\theta\frac{d}{dt} \int_{\mathbb{R}^2}  \phi^2 \text{div}u \text{d}x-\theta\int_{\mathbb{R}^2}  (\phi^2)_t \text{div}u \text{d}x\\
=&\theta\frac{d}{dt} \int_{\mathbb{R}^2}  \phi^2 \text{div}u \text{d}x+\theta\int_{\mathbb{R}^2}  u \cdot \nabla \phi^2 \text{div}u \text{d}x+\theta(\gamma-1)\int_{\mathbb{R}^2}   \phi^2 (\text{div}u)^2 \text{d}x\\
\leq &\theta\frac{d}{dt} \int_{\mathbb{R}^2}  \phi^2 \text{div}u \text{d}x+C(|\nabla u|^2_2+|\text{div}u|_{\infty}+|\text{div}u|_{\infty}|\nabla\phi|_2^{2}),\\
L_{9}=&-\int_{\mathbb{R}^2} (\psi\cdot Q(u)) \cdot (\alpha \triangle u+(\alpha+\beta)\nabla \text{div}u)\text{d}x\\
\leq& C|\psi|_6|\nabla^2 u|^{\frac{4}{3}}_2|\nabla u|^{\frac{2}{3}}_2\leq C(\epsilon)|\nabla u|^2_2+\epsilon |\nabla^2 u|^2_2,
\end{split}
\end{equation}
where $\epsilon> 0$ is a sufficiently small constant.
Combining (\ref{zhu6})-(\ref{zhu10t33}),  we have
\begin{equation}\label{zhu6qss}
\begin{split}
&\frac{1}{2} \frac{d}{dt}\int_{\mathbb{R}^2}(\alpha|\nabla u|^2+(\alpha+\beta)|\text{div}u|^2-\theta \phi^2 \text{div}u\Big)\text{d}x+C|\nabla^2 u|^2_2\\
\leq &C((|\nabla u|^2_2+|\nabla \phi|^2_2)(|\text{div}u|_\infty+1)+|\text{div}u|_\infty).
\end{split}
\end{equation}

Second, applying $\nabla$ to  $(\ref{we22})_1$ and multiplying by $(\nabla \phi)^{\top}$,  we have
\begin{equation}\label{zhu20}
\begin{split}
&(|\nabla \phi|^2)_t+\text{div}(|\nabla \phi|^2u)+(\gamma-2)|\nabla \phi|^2\text{div}u\\
=&-2 (\nabla \phi)^\top \nabla u( \nabla \phi)-(\gamma-1) \phi \nabla \phi \cdot \nabla \text{div}u\\
=&-2 (\nabla \phi)^\top D(u) (\nabla \phi)-(\gamma-1) \phi \nabla \phi \cdot \nabla \text{div}u.
\end{split}
\end{equation}
Integrating (\ref{zhu20}) over $\mathbb{R}^2$, we get
\begin{equation}\label{zhu21}
\begin{split}
\frac{d}{dt}|\nabla \phi|^2_2
\leq& C(\epsilon)(|D( u)|_\infty+1)|\nabla \phi|^2_2+\epsilon |\nabla^2 u|^2_2.
\end{split}
\end{equation}
Adding (\ref{zhu21}) to (\ref{zhu6qss}), from Gronwall's inequality we immediately obtain
\begin{equation*}
\begin{split}
|\nabla u(t)|^2_{ 2}+|\nabla \phi(t)|^2_{ 2}+\int_0^t |\nabla^2 u(s)|^2_2\text{d}t\leq C,\quad 0\leq t\leq T.
\end{split}
\end{equation*}

Finally, due to $u_t=-Lu-u\cdot \nabla u-2\theta\phi\nabla \phi+\psi\cdot Q(u)$, we deduce that
\begin{equation*}
\begin{split}
\int_0^t | u_t|^2_2\text{d}t\leq C\int_0^t (| \nabla^2 u|^2_2+|\nabla u|^2_3|u|^2_6+|\phi|^2_\infty|\nabla \phi|^2_2+|\nabla u|^2_3|\psi|^2_6)\text{d}t\leq C.
\end{split}
\end{equation*}
 \end{proof}

Next, we proceed to improve the regularity of $\phi$, $\psi$ and $u$. To this end, we first drive some bounds on derivatives of $u$ based on the above estimates.
 \begin{lemma}\label{s6} Let $( \rho, u)$ be the unique regular solution to  the Cauchy problem  (\ref{eq:1.1}) on $[0,\overline{T})$ satisfying \eqref{we11}. Then
 \begin{equation}\label{zhu14ss}
\begin{split}
&\sup_{0\leq t\leq T}|u_t(t)|^2_2+\sup_{0\leq t\leq T}|u(t)|_{D^2}+\int_0^T|\nabla u_t|^2_2\text{d}t\leq C,\quad 0\leq T<  \overline{T},
\end{split}
\end{equation}
where $C$ only depends on $C_0$ and $T$.
\end{lemma}
\begin{proof}
Using $Lu=-u_t-u\cdot \nabla u-2\theta\phi\nabla \phi+\psi\cdot Q(u)$ and Lemma \ref{zhenok}, we have
 \begin{equation}\label{zhu15nn}
\begin{split}
|u|_{D^2}\leq& C(|u_t|_2+|u|_6|\nabla u|^{\frac{2}{3}}_2|\nabla^2 u|^{\frac{1}{3}}+|\phi|_\infty|\nabla \phi|_2+ |\psi|_6|\nabla u|^{\frac{2}{3}}_2|\nabla^2 u|^{\frac{1}{3}}_2),
\end{split}
\end{equation}
which implies, from Young's inequality with appropriate weights, that
\begin{equation}\label{zhu15nn1}
\begin{split}
|u|_{D^2}\leq& C(|u_t|_2+|u|^{\frac{3}{2}}_6|\nabla u|_2+|\nabla \phi|_2+ |\nabla u|_2)\leq C(1+|u_t|_2).
\end{split}
\end{equation}

Next, differentiating $(\ref{we22})_3$ with respect to $t$, it reads
\begin{equation}\label{zhu37ss}
\begin{split}
u_{tt}+Lu_t=-(u\cdot\nabla u)_t -2\theta(\phi\nabla \phi)_t+(\psi\cdot Q(u))_t.
\end{split}
\end{equation}
Multiplying (\ref{zhu37ss}) by $u_t$ and integrating over $\mathbb{R}^2$, one has
\begin{equation}\label{zhu38}
\begin{split}
&\frac{1}{2} \frac{d}{dt}|u_t|^2_2+\alpha|\nabla u_t|^2_2+(\alpha+\beta)|\text{div} u_t|^2_2\\
=&\int_{\mathbb{R}^2} \Big(-(u\cdot \nabla u)_t \cdot u_t-2\theta(\phi \nabla \phi)_t \cdot u_t+(\psi \cdot Q(u))_t\cdot u_t \Big) \text{d}x\equiv: \sum_{i=10}^{12} L_i.
\end{split}
\end{equation}
We estimate the right-hand side of \eqref{zhu38} term by term as follows.
\begin{equation}\label{zhu9fg}
\begin{split}
L_{10}=&-\int_{\mathbb{R}^2} (u\cdot \nabla u)_t \cdot u_t \text{d}x=-\int_{\mathbb{R}^2} \Big(\big(u_t \cdot \nabla u\big) \cdot u_t +\big(u \cdot \nabla u_t \big)\cdot u_t\Big)\text{d}x\\
=&-\int_{\mathbb{R}^2} \big(u_t \cdot D( u) \cdot u_t -( u_t)^2 \text{div}u\big)\text{d}x\leq C|D(u)|_\infty| u_t|^2_2,\\
L_{11}=&-\int_{\mathbb{R}^2} 2\theta(\phi \nabla \phi)_t \cdot u_t \text{d}x=\theta\int_{\mathbb{R}^2}   (\phi^2)_t \text{div}u_t \text{d}x\\
=&-\frac{\theta(\gamma-1)}{2}\frac{d}{dt}\int_{\mathbb{R}^2}   \phi^2 (\text{div}u)^2 \text{d}x-\frac{\theta(\gamma-1)}{2}\int_{\mathbb{R}^2}   u\phi^2 \nabla(\text{div}u)^2 \text{d}x\\
&-\frac{\theta(\gamma-1)^2}{2}\int_{\mathbb{R}^2}   \phi^2 (\text{div}u)^3 \text{d}x-\theta\int_{\mathbb{R}^2}   u\cdot \nabla \phi^2 \text{div}u_t \text{d}x\\
\leq&-\frac{\theta(\gamma-1)}{2}\frac{d}{dt}\int_{\mathbb{R}^2}   \phi^2 (\text{div}u)^2 \text{d}x+ C|u|_\infty| \phi|^2_\infty|\nabla u|_2|\nabla^2 u|_2\\
&+C|\phi|^2_\infty|D( u)|_\infty|\nabla u|^2_2+C|\phi|_\infty|\nabla \phi|_2|u|_\infty|\nabla u_t|_2 \\
\leq&-\frac{\theta(\gamma-1)}{2}\frac{d}{dt}\int_{\mathbb{R}^2}   \phi^2 (\text{div}u)^2 \text{d}x+ \frac{\alpha}{4}|\nabla u_t|^2_2+ C(1+|D( u)|_\infty+|u|^2_{D^2}),\\
L_{12}=&\int_{\mathbb{R}^2} (\psi \cdot Q(u))_t\cdot u_t\text{d}x=\int_{\mathbb{R}^2} \psi \cdot Q(u)_t\cdot u_t\text{d}x\\
&-\int_{\mathbb{R}^2} \nabla \text{div}u \cdot Q(u)\cdot u_t\text{d}x+\int_{\mathbb{R}^2} \psi\cdot u \text{div}(Q(u)\cdot u_t)\text{d}x\\
\leq& C|\psi|_6|\nabla u_t|_2|u_t|^{\frac{2}{3}}_2| \nabla u_t|^{\frac{1}{3}}_2+C|\nabla^2 u|_2 |Q( u)|_\infty |u_t|_2\\
&+C|\psi|_6| u|_6|\nabla^2 u|_2  | u_t|_6+C|\psi|_6| u|_6||Q(u)|_6  |\nabla u_t|_2\\
\leq&\frac{\alpha}{10}|\nabla u_t|^2_2+C(1+ |D( u)|_\infty)(|u_t|^2_2+|u|^2_{D^2}).
\end{split}
\end{equation}
It is clear from (\ref{zhu38})-(\ref{zhu9fg}), \eqref{we11}, and \eqref{zhu15nn1} that
\begin{equation}\label{zhu8qqss}
\begin{split}
& \frac{d}{dt}(|u_t|^2_2+|\phi \text{div}u|^2_2)+|\nabla u_t|^2_2\\
\leq &C(1+ |D( u)|_\infty)|u_t|^2_2+C(1+|D( u)|_\infty).
\end{split}
\end{equation}
Integrating (\ref{zhu8qqss}) over $(\tau,t)$ $(\tau \in( 0,t))$ , we have
\begin{equation}\label{zhu13vbn}
\begin{split}
&|u_t(t)|^2_2+|\phi \text{div}u(t)|^2_2+\int_\tau^t|\nabla u_t(s)|^2_2\text{d}s\\
\leq & |u_t(\tau)|^2_2+|\phi \text{div}u(\tau)|^2_2+C\int_{\tau}^t  \Big((|D( u)|_\infty+1)|u_t|^2_2\Big)(s)\text{d}s+C.
\end{split}
\end{equation}

From the momentum equations $(\ref{we22})_3$, we obtain
\begin{equation}\label{zhu15vbnm}
\begin{split}
|u_t(\tau)|_2\leq C\big( |u|_\infty |\nabla u|_2+|\phi|_\infty|\nabla \phi|_2+|u|_{D^2}+|\psi|_6|\nabla u|_2\big)(\tau),
\end{split}
\end{equation}
which, together with (\ref{reg11}), gives
\begin{equation}\label{zhu15vbnm}
\begin{split}
\lim \sup_{\tau\rightarrow 0}|u_t(\tau)|_2\leq C\big( |u_0|_\infty |\nabla u_0|_2+|\phi_0|_\infty|\nabla \phi_0|_2+|u_0|_{D^2}+|\psi_0|_6|\nabla u_0|_2\big)\leq C_0.
\end{split}
\end{equation}
Letting $\tau \rightarrow 0$ in (\ref{zhu13vbn}),  applying Gronwall's inequality, we arrive at
$$
|u_t(t)|^2_2+\int_0^t|\nabla u_t(s)|^2_2\text{d}s+|u(t)|_{D^2}\leq C, \quad 0\leq t\leq T.
$$
This completes the proof of this lemma.
\end{proof}

The following lemma gives bounds of $\nabla \phi$ and $\nabla^2 u$.
 \begin{lemma}\label{s7} Let $( \rho, u)$ be the unique regular solution to  the Cauchy problem  (\ref{eq:1.1}) on $[0,\overline{T})$ satisfying \eqref{we11}. Then
 \begin{equation}\label{zhu54}
\begin{split}
&\sup_{0\leq t\leq T}\|\phi(t)\|_{W^{1,6}}+\sup_{0\leq t\leq T}|\phi_t(t)|_6+\int_0^T|u(t)|^2_{D^{2,6}}\text{d}t\leq C,
\quad 0\leq T<  \overline{T},
\end{split}
\end{equation}
where $C$ only depends on $C_0$ and $T$.
\end{lemma}
 \begin{proof}
 First, using Lemma \ref{zhenok} , we read from $Lu=-u_t-u\cdot \nabla u-2\theta\phi\nabla \phi+\psi\cdot Q(u)$ that
 \begin{equation}\label{zhu55}
\begin{split}
|\nabla^2 u|_6 \leq& C(|u_t|_6+|u\cdot \nabla u|_6+|\phi\nabla \phi|_6+|\psi \cdot Q( u)|_6)\\
\leq& C(1+|\nabla u_t|_2+|\nabla \phi|_6+ |D( u)|^{\frac{2}{5}}_2 |\nabla D(u)|^{\frac{3}{5}}_6),
\end{split}
\end{equation}
where we have used the fact that $|D( u)|_\infty\leq C|D( u)|^{\frac{2}{5}}_2 |\nabla D(u)|^{\frac{3}{5}}_6$. Now, Young's inequality implies that
 \begin{equation}\label{zhu65sss}
\begin{split}
|\nabla^2 u|_6
\leq C(1+|\nabla u_t|_2+|\nabla \phi|_6).
\end{split}
\end{equation}

Next,  applying $\nabla$ to  $(\ref{we22})_1$, multiplying the result equations by $6|\nabla \phi|^4 \nabla \phi$, we have
\begin{equation}\label{zhu20cccc}
\begin{split}
&(|\nabla \phi|^6)_t+\text{div}(|\nabla \phi|^6u)+(3\gamma-4)|\nabla \phi|^6\text{div}u\\
=&-6 |\nabla \phi|^4(\nabla \phi)^\top \nabla u (\nabla \phi)-(3\gamma-3) \phi|\nabla \phi|^4 \nabla \phi \cdot \nabla \text{div}u\\
=&-6 |\nabla \phi|^4(\nabla \phi)^\top D( u) (\nabla \phi)-(3\gamma-3) \phi|\nabla \phi|^4 \nabla \phi \cdot \nabla \text{div}u,
\end{split}
\end{equation}
which implies, upon integrating over $\mathbb{R}^2$, that
\begin{equation}\label{zhu200}
\begin{split}
\frac{d}{dt}|\nabla \phi|_6
\leq& C(|D( u)|_\infty)|\nabla \phi|_6+C|\nabla^2 u|_6.
\end{split}
\end{equation}
Therefore, we obtain from \eqref{zhu65sss} that
\begin{equation}\label{zhu28ss}
\begin{split}
\frac{d}{dt}|\nabla \phi|_6
\leq& C(1+|D( u)|_\infty)|\nabla \phi|_6+C(1+|\nabla u_t|^2_2).
\end{split}
\end{equation}
In view of Lemma \ref{s6} and \eqref{we11} , we apply Gronwall's inequality to conclude that
$$
|\nabla \phi(t)|_6\leq C(1+|\nabla \phi_0|_6)\exp \Big(\int_0^t(1+|D( u)|_\infty) \text{d}s\Big)\leq C, \quad 0\leq t \leq T.
$$

Finally,  we infer from (\ref{zhu65sss}) and Lemma \ref{s6} that
\begin{equation}\label{zhu25ss}
\begin{split}
\int_0^t|u(s)|^2_{D^{2,6}}\text{d}s
\leq& C\int_0^t(1+|\nabla \phi(s)|^2_6+ |\nabla u_t(s)|^2_2)\text{d}s\leq C,\quad 0\leq t\leq T.
\end{split}
\end{equation}
The proof of the lemma is finished.
 \end{proof}

Lemma \ref{s7} implies that
 \begin{equation}\label{keygradient}
\int_0^t|\nabla u(\cdot, s)|_\infty \text{d}s\leq C, \end{equation}
for any $t\in [0, \overline{T})$ with $C>0$ a finite number. Noting that \eqref{we22} is essentially a parabolic-hyperbolic system, it is then standard to derive other higher order estimates for the regularity of
the regular solutions. We will show this fact in the following 4 lemmas.

\begin{lemma}\label{s8}  Let $( \rho, u)$ be the unique regular solution to  the Cauchy problem  (\ref{eq:1.1}) on $[0,\overline{T})$ satisfying \eqref{we11}. Then
\begin{equation*}
\begin{split}
&\sup_{0\leq t\leq T}|\phi(t)|^2_{D^2}+\sup_{0\leq t\leq T}|\psi(t)|^2_{D^1}+\sup_{0\leq t\leq T}\|\phi_t(t)\|^2_1+\sup_{0\leq t\leq T}|\psi_t(t)|^2_{2}\\
&\quad+\int_{0}^{T}\Big(|u(t)|^2_{D^{3}}+|\phi_{tt}(t)|^2_{2}\Big)\text{d}t\leq C, \quad 0\leq T<  \overline{T},
\end{split}
\end{equation*}
where $C$ only depends on $C_0$ and $T$.
 \end{lemma}
 \begin{proof} Using $Lu=-u_t-u\cdot \nabla u-2\theta\phi\nabla \phi+\psi\cdot Q(u)$ and Lemma \ref{zhenok}, we have
 \begin{equation}\label{zhu150}
\begin{split}
|u|_{D^3}\leq& C(|u_t|_{D^1}+|u\cdot \nabla u|_{D^1}+|\phi\nabla \phi|_{D^1}+|\psi \cdot Q( u)|_{D^1})\\
\leq& C(1+|u_t|_{D^1}+| \phi|_{D^2}+|\psi|_6 |\nabla^2 u|_3+|\nabla\psi|_2 |D( u)|_\infty)\\
\leq& C(1+|u_t|_{D^1}+| \phi|_{D^2}+|\psi|_6 |\nabla^2 u|^{\frac{2}{3}}_2|\nabla^3 u|^{\frac{1}{3}}_2+|\nabla\psi|_2 |D( u)|^{\frac{2}{3}}_6|\nabla D( u)|^{\frac{1}{3}}_6),
\end{split}
\end{equation}
where we have used the fact that $|D( u)|_\infty\leq C |D( u)|^{\frac{2}{3}}_6|\nabla D( u)|^{\frac{1}{3}}_6$. With the help of Young's inequality, (\ref{zhu150}) offers that
 \begin{equation}\label{zhu1500}
|u|_{D^3}
\leq C(1+|u_t|_{D^1}+| \phi|_{D^2}+|D( u)|_\infty|\nabla\psi|_2).
\end{equation}

Next, applying $\partial_i$ ($i=1,2$) to  $(\ref{we22})_2$  with respect to $x$, we obtain
\begin{equation}\label{hypss}\begin{split}
&(\partial_i \psi)_t+\sum_{l=1}^2 A_l \partial_l\partial_i \psi+B\partial_i \psi+\partial_i \nabla \text{div} u \\
=&\Big(-\partial_i(B\psi)+B\partial_i \psi\Big)+\sum_{l=1}^2 \Big(-\partial_i(A_l) \partial_l \psi\Big)=\Theta'_1+\Theta'_2.
\end{split}
\end{equation}
Multiplying (\ref{hypss}) by $2(\partial_i\psi)^{\top}$, integrating over $\mathbb{R}^2$, and then summing over $i$, noting that  $A_l$ ($l=1,2$) are symmetric, it is not difficult to show that
\begin{equation}\label{zhenzhen11}\begin{split}
\frac{d}{dt}|\nabla \psi|^2_2\leq& C\big|\text{div}A\big|_\infty|\nabla \psi|^2_2+\int_{\mathbb{R}^2}\sum_{i=1}^2(\partial_i\psi)^\top (\nabla u)^\top( \partial_i \psi)+C|\nabla^3 u|_2|\nabla \psi|_2\\
&+C|\Theta'_1 |_2|\nabla\psi|_2+2\int_{\mathbb{R}^2}\sum_{i=1}^2 (\partial_i\psi)^{\top}\Theta'_2\text{d}x,
\end{split}
\end{equation}
where $\text{div}A=\displaystyle\sum_{l=1}^{2}\partial_{l}A_l$. We treat each term on the right-hand side of the above inequality as follows.  From the definition of matrices $A_l$ and $B$, it is clear that
\begin{equation}
\big|\text{div}A\big|_\infty|\nabla \psi|^2_2+\int_{\mathbb{R}^2}\sum_{i=1}^2(\partial_i\psi)^\top (\nabla u)^\top (\partial_i \psi)
\leq  C|\nabla u|_\infty |\nabla \psi|^2_2.
\end{equation}
When $|\zeta|=1$, choosing $r=2,\ a=3$, $b=6$ in (\ref{ku22}), we have
\begin{equation}\label{zhen211}
|\Theta'_1|_2=|D^\zeta(B\psi)-BD^\zeta \psi|_2\leq C|\nabla^2 u|_3|\psi|_6.
\end{equation}
For the last term on the right-hand side of (\ref{zhenzhen11}), we have
\begin{equation}\label{ghcc}
2\int_{\mathbb{R}^2}\sum_{i=1}^2 (\partial_i\psi)^{\top}\Theta'_2\text{d}x\le  C|\nabla u|_\infty|\nabla \psi|^2_2.
\end{equation}
Combining (\ref{zhu1500}), \eqref{zhenzhen11}-(\ref{ghcc}), and using Gagliardo-Nirenberg inequality, we have
\begin{equation}\label{pl}
\begin{split}
\frac{d}{dt}|\psi|^2_{D^1}\leq& C(1+|\nabla u|_\infty)|\psi|^2_{D^1}+C(1+|\nabla^3 u|_2)|\psi|_{D^1}+C\\
\leq& C(1+|\nabla u|_\infty)|\psi|^2_{D^1}+C(1+|\phi|^2_{D^2}+|\nabla u_t|^2_2).
\end{split}
\end{equation}

On the other hand, let $\nabla \phi=G=(G^{(1)},G^{(2)})^\top$. Applying $\nabla^2$ to $(\ref{we22})_1$, we have
\begin{equation}\label{zhuu11}
\begin{split}
0=&(\nabla G)_t+\nabla((\nabla u)^\top \cdot G)+\nabla ((\nabla G)^\top\cdot u)+\frac{\gamma-1}{2}\nabla(G\text{div}u+\phi \nabla \text{div}u\big)\\
=&(\nabla G)_t+\sum_{i=1}^2(G^{(i)}\nabla^2 u^{(i)}+\nabla u^{(i)}\otimes \nabla G^{(i)})+\sum_{i=1}^2(u^{(i)}\nabla^2G^{(i)}+\nabla G^{(i)}\otimes \nabla u^{(i)})\\
&+\frac{\gamma-1}{2}\big(\nabla G\text{div}u+ G \otimes\nabla\text{div}u\big)+\frac{\gamma-1}{2}\big(\phi \nabla^2 \text{div}u+\nabla \text{div}u\otimes G\big).
\end{split}
\end{equation}
Similarly to the previous step, we multiply (\ref{zhuu11}) by $ 2\nabla G$ and integrate it over ${\mathbb {R}}^2$ to derive

\begin{equation}\label{zhuu12acc}
\begin{split}
\frac{d}{dt}|G|^2_{D^1}
\leq& C|\nabla u|_\infty|G|^2_{D^1}+C|G|_6|\nabla G|_2|\nabla^2 u|_3+C|\phi|_\infty|\nabla G|_2|\nabla^3 u|_2\\
\leq& C(1+|\nabla u|_\infty)(|G|^2_{D^1}+|\psi|^2_{D^1})+C(1+|\nabla u_t|^2_2).
\end{split}
\end{equation}
This estimate, together with \eqref{pl}, gives that
\begin{equation}\label{zhmm12acc}
\begin{split}
\frac{d}{dt}(|G|^2_{D^1}+|\psi|^2_{D^1})
\leq& C(1+|\nabla u|_\infty)(|G|^2_{D^1}+|\psi|^2_{D^1})+C(1+|\nabla u_t|^2_2).\end{split}
\end{equation}
Then the Gronwall's inequality and \eqref{keygradient} imply
\begin{equation*}
\begin{split}
|\phi(t)|^2_{D^2}+|\psi(t)|^2_{D^1}+\int_{0}^{t}|u(s)|^2_{D^{3}}\text{d}t\leq C, \quad 0\leq t\leq T.
\end{split}
\end{equation*}

Finally, using the following relations
\begin{equation}\label{ghtuss}
\begin{split}
\psi_t=&-\nabla (u \cdot \psi)-\nabla \text{div} u,\
\phi_t=-u\cdot \nabla \phi-\frac{\gamma-1}{2}\phi\text{div} u,\\
\phi_{tt}=&-u_t\cdot \nabla \phi-u\cdot \nabla \phi_t-\frac{\gamma-1}{2}\phi_t\text{div} u-\frac{\gamma-1}{2}\phi\text{div} u_t,
\end{split}
\end{equation}
we conclude the proof of this lemma.
\end{proof}
In order to obtain higher order regularity, we need the following improved estimate.
\begin{lemma}\label{s9} Let $( \rho, u)$ be the unique regular solution to  the Cauchy problem  (\ref{eq:1.1}) on $[0,\overline{T})$ satisfying \eqref{we11}. Then
\begin{equation*}
\begin{split}
\sup_{0\leq t\leq T}|u_t(t)|^2_{D^1}+\sup_{0\leq t\leq T}|u(t)|^2_{D^3}+\int_{0}^{T}|u_{tt}(t)|^2_2\text{d}t\leq C, \quad 0\leq T<  \overline{T},
\end{split}
\end{equation*}
where $C$ only depends on $C_0$ and $T$.
 \end{lemma}
\begin{proof} First, $Lu_t=-u_{tt}-(u\cdot \nabla u)_t-2\theta(\phi\nabla \phi)_t+(\psi\cdot Q(u))_t$ and Lemma \ref{zhenok} yield
 \begin{equation}\label{zhu150z}
\begin{split}
|u_t|_{D^2}\leq& C(|u_{tt}|_{2}+|(u\cdot \nabla u)_t|_{2}+|\nabla (\phi^2)_t|_{2}+|(\psi \cdot Q( u))_t|_{2})\\
\leq& C(|u_{tt}|_{2}+|u|_{\infty}|\nabla u_t|_2+|\nabla u|_3 | u_t|^{\frac{1}{3}}_2 |\nabla u_t|^{\frac{2}{3}}_2+|\phi|_\infty|\nabla\phi_t|_2)\\
&+C(|\nabla\phi|_3|\phi_t|_6+|\psi|_6|\nabla u_t|^{\frac{2}{3}}_2|\nabla^2 u_t|^{\frac{1}{3}}_2+|\psi_t|_2|Q( u)|_\infty)\\
\leq& C(1+|u_{tt}|_{2}+|\nabla u_t|_2+|u|_{D^{2,6}}),
\end{split}
\end{equation}
which implies, with the help of Young's inequality, that
 \begin{equation}\label{zhu1500x}
|u_t|_{D^2}
\leq C(1+|u_{tt}|_{2}+|\nabla u_t|_2+|u|_{D^{2,6}}).
\end{equation}

Now, multiplying  $(\ref{zhu37ss})$ by $u_{tt}$ and integrating over $\mathbb{R}^2$,  we have
\begin{equation}\label{zhu19ss}
\begin{split}
& \frac{1}{2}\frac{d}{dt}\Big(\alpha|\nabla u_t|^2_2+(\alpha+\beta)|\text{div}u_t|^2_2\Big)+| u_{tt}|^2_2\\
=&\int_{\mathbb{R}^2} \Big(-(u\cdot \nabla u)_t\cdot u_{tt}-2\theta(\phi \nabla \phi)_t \cdot u_{tt}+(\psi \cdot Q(u))_t\cdot u_{tt} \Big) \text{d}x\equiv: \sum_{i=13}^{15} L_i.
\end{split}
\end{equation}
For the terms $L_{13}$--$L_{15}$, we perform the following estimates:
\begin{equation}\label{zhu201}
\begin{split}
L_{13}=&-\int_{\mathbb{R}^2} (u\cdot \nabla u)_t \cdot u_{tt} \text{d}x\leq C|u_t|_{6}|\nabla u|_{3}|u_{tt}|_2+C|u|_\infty |\nabla u_t|_2|u_{tt}|_2\\
\leq& C(|\nabla u|^2_{3}+|u|^2_\infty )|\nabla u_t|^2_2 +\frac{1}{10}|u_{tt}|^2_2,\\
L_{14}=&-\int_{\mathbb{R}^2} 2\theta(\phi \nabla \phi)_t \cdot u_{tt} \text{d}x
=\theta\frac{d}{dt}\int_{\mathbb{R}^2} (\phi^2)_t \cdot \text{div}u_{t} \text{d}x-\int_{\mathbb{R}^2} \theta(\phi^2)_{tt} \cdot \text{div}u_{t} \text{d}x\\
\leq&\theta\frac{d}{dt}\int_{\mathbb{R}^2} (\phi^2)_t \cdot \text{div}u_{t} \text{d}x+ C|(\phi^2)_{tt}|^2_2+C |\nabla u_t|^2_2,\\
L_{15}=&\int_{\mathbb{R}^2} (\psi \cdot Q(u))_t\cdot u_{tt}\text{d}x
=\int_{\mathbb{R}^2} \psi \cdot Q(u)_t\cdot u_{tt}\text{d}x+\int_{\mathbb{R}^2} \psi_t \cdot Q(u)\cdot u_{tt}\text{d}x\\
\leq& C|\psi|_6 |\nabla u_t|_3 |u_{tt}|_2+|\psi_t|_2|\nabla u|_\infty|u_{tt}|_2\\
\leq& C|\psi|^2_6|\nabla u_t|^{\frac{4}{3}}_2|\nabla^2 u_t|^{\frac{2}{3}}_2 +C|\psi_t|^2_2 |\nabla u|^{\frac{4}{5}}_2|\nabla^2 u|^{\frac{6}{5}}_6  +\frac{1}{10}|u_{tt}|^2_2\\
\leq& C(1+|\nabla u_t|^2_2+|\nabla^2 u|^{2}_6) +\frac{1}{5}|u_{tt}|^2_2.
\end{split}
\end{equation}
Therefore, \eqref{zhu19ss} and \eqref{zhu201} imply that
\begin{equation}\label{zhu19qqss}
\begin{split}
& \frac{1}{2}\frac{d}{dt}\Big(|\nabla u_t|^2_2+(\alpha+\beta)|\text{div}u_t|^2_2-\theta\int_{\mathbb{R}^2} (\phi^2)_t \cdot \text{div}u_{t} \text{d}x\Big)+| u_{tt}|^2_2\\
\leq& C(1+|\nabla u_t|^2_2+|\nabla^2 u|^{2}_6),
\end{split}
\end{equation}
which, upon integrating over $(\tau,t)$, yields
\begin{equation}\label{liu03}
\begin{split}
|\nabla u_t(t)|^2_2+\int_{\tau}^t| u_{tt}(s)|^2_2\text{d}s\leq C +|\nabla u_t(\tau)|^2_2+\int_{\tau}^t| \nabla u_{t}|^2_2\text{d}s, \quad  0\leq t \leq T,
\end{split}
\end{equation}
where we used the fact that for any $\epsilon>0$
\begin{equation}\label{liu01}
\begin{split}
\int_{\mathbb{R}^3} (\phi^2)_t \cdot \text{div}u_{t} \text{d}x
\leq& \epsilon|\nabla u_t|^2_2+C.
\end{split}
\end{equation}
From the momentum equations $ (\ref{we22})_3$, we  have
\begin{equation}\label{liu01}
\begin{split}
|\nabla u_t(\tau)|_2\leq C\big(1+ \|u\|^2_3+\|\phi\|^2_3+|\psi|_{L^6\cap D^1\cap D^2}\| u\|_2\big)(\tau).
\end{split}
\end{equation}
Using the regularity in (\ref{reg11}), we find
\begin{equation}\label{liu02}
\begin{split}
\lim \sup_{\tau\rightarrow 0}|\nabla u_t(\tau)|_2\leq C\big(1+ \|u_0\|^2_3+\|\phi_0\|^2_3+|\psi_0|_{L^6\cap D^1\cap D^2}\| u_0\|_2\big)\leq C_0.
\end{split}
\end{equation}
Letting $\tau \rightarrow 0$ in (\ref{liu03}), we finally proved that
\begin{equation}\label{zhu22wsx}
\begin{split}
|\nabla u_t(t)|^2_2+\int_{0}^t| u_{tt}(s)|^2_2\text{d}s\leq C, \quad  0\leq t \leq T.
\end{split}
\end{equation}
In order to complete the proof of this lemma, we observe from  (\ref{zhu150}) and (\ref{zhu1500x}) that
\begin{equation*}
\begin{split}
&|u(t)|_{D^3}\leq C \big(1+|u_t|_{D^1}\big)(t)\leq C,\quad  0\leq t \leq T.\\
&\int_0^t |u_t|^2_{D^2}\text{d}s\leq \int_0^t C\big(1+|u_{tt}|^2_{2}+|\nabla u_t|^2_2+|u|^2_{D^{2,6}}\big)\text{d}s\leq C, \quad  0\leq t \leq T.
\end{split}
\end{equation*}
\end{proof}

It remains to prove the following lemma for the required regularity estimate.

\begin{lemma}\label{s10} Let $( \rho, u)$ be the unique regular solution to  the Cauchy problem  (\ref{eq:1.1}) on $[0,\overline{T})$ satisfying \eqref{we11}. Then
\begin{equation*}
\begin{split}
&\sup_{0\leq t\leq T}(|\phi(t)|^2_{D^3}+|\psi(t)|^2_{D^2}+\|\phi_t(t)\|^2_2+|\psi_t(t)|^2_{D^1})+\int_{0}^{T}|u|^2_{D^{4}}\text{d}t\leq C,\\
&\sup_{0\leq t\leq T}(t|u_t(t)|^2_{ D^2}+t|u_{tt}(t)|^2_2+t|u(t)|^2_{D^4})+\int_{0}^{T}(t|u_{tt}|^2_{D^1}+t|u_{t}|^2_{D^3})\text{d}t\leq C,
\end{split}
\end{equation*}
where $0\leq T<  \overline{T}$,  and  $C$ only depends on $C_0$ and $T$.
 \end{lemma}
 \begin{proof} The first assertion in this lemma follows in the similar lines of proof as that for Lemma \ref{s7}, while the second assertion can be proved by the same method used in Lemma \ref{5} and Lemma \ref{s9}. We omit the details for
simplicity of presentation.
\end{proof}

Now we know from Lemmas \ref{s1}-\ref{s10} that, if the regular solution $(\rho, u)(x,t)$ exists up to the time ${\overline{T}}>0$, with the maximal time ${\overline{T}}<+\infty$ such that  the assumption \eqref{we11} holds,
then $(\rho^{\frac{\gamma-1}{2}},\nabla \rho/\rho,u)|_{t=\overline{T}} =\lim_{t\rightarrow \overline{T}}(\rho^{\frac{\gamma-1}{2}},\nabla \rho/\rho,u)$ satisfy the conditions imposed on the initial data $(\ref{th78})$. If we solve the system (\ref{eq:1.1}) with the initial time ${\overline{T}}$, then Theorem 1.1 ensures that $(\rho, u)(x,t)$ extends beyond ${\overline{T}}$ as the unique regular solution. This contradicts to the fact that ${\overline{T}}$ is the maximal existence time. We thus complete the proof of Theorem \ref{th3s}.

\bigskip

{\bf Acknowledgement:} The research of  Y. Li and S. Zhu were supported in part
by National Natural Science Foundation of China under grant 11231006 and Natural Science Foundation of Shanghai under grant 14ZR1423100. S. Zhu was also supported by China Scholarship Council. The research of R. Pan was partially supported
by National Science Foundation under grants DMS-0807406 and DMS-1108994.

\end{document}